# RENEWAL THEORY AND COMPUTABLE CONVERGENCE RATES FOR GEOMETRICALLY ERGODIC MARKOV CHAINS

By Peter H. Baxendale

*University of Southern California*

We give computable bounds on the rate of convergence of the transition probabilities to the stationary distribution for a certain class of geometrically ergodic Markov chains. Our results are different from earlier estimates of Meyn and Tweedie, and from estimates using coupling, although we start from essentially the same assumptions of a drift condition toward a "small set." The estimates show a noticeable improvement on existing results if the Markov chain is reversible with respect to its stationary distribution, and especially so if the chain is also positive. The method of proof uses the first-entrance–last-exit decomposition, together with new quantitative versions of a result of Kendall from discrete renewal theory.

**1. Introduction.** Let $\{X_n : n \geq 0\}$ be a time homogeneous Markov chain on a state space $(S, \mathcal{B})$. Let $P(x, A), x \in S, A \in \mathcal{B}$ denote the transition probability and let $P$ denote the corresponding operator on measurable functions $S \to \mathbf{R}$. There has been much interest and activity recently in obtaining computable bounds for the rate of convergence of the time $n$ transition probability $P_n(x, \cdot)$ to a (unique) invariant probability measure $\pi$. These estimates are of importance for simulation techniques such as Markov chain Monte Carlo (MCMC).

Throughout this paper we assume the following conditions are satisfied.

(A1) *Minorization condition.* There exist $C \in \mathcal{B}$, $\tilde{\beta} > 0$ and a probability measure $\nu$ on $(S, \mathcal{B})$ such that

$$P(x, A) \geq \tilde{\beta} \nu(A)$$

for all $x \in C$ and $A \in \mathcal{B}$.

---

(A2) *Drift condition.* There exist a measurable function $V\colon S \to [1, \infty)$ and constants $\lambda < 1$ and $K < \infty$ satisfying

$$PV(x) \le \begin{cases} \lambda V(x), & \text{if } x \notin C, \\ K, & \text{if } x \in C. \end{cases}$$

(A3) *Strong aperiodicity condition.* There exists $\beta > 0$ such that $\tilde{\beta}\nu(C) \ge \beta$.

The following result converts information about the one-step behavior of the Markov chain into information about the long term behavior of the chain.

THEOREM 1.1. *Assume* (A1)–(A3). *Then* $\{X_n\colon n \ge 0\}$ *has a unique stationary probability measure* $\pi$, *say, and* $\int V \, d\pi < \infty$. *Moreover, there exists* $\rho < 1$ *depending only (and explicitly) on* $\beta$, $\tilde{\beta}$, $\lambda$ *and* $K$ *such that whenever* $\rho < \gamma < 1$ *there exists* $M < \infty$ *depending only (and explicitly) on* $\gamma$, $\beta$, $\tilde{\beta}$, $\lambda$ *and* $K$ *such that*

$$(1) \qquad \sup_{|g| \le V} \left| (P^n g)(x) - \int g \, d\pi \right| \le MV(x)\gamma^n$$

*for all* $x \in S$ *and* $n \ge 0$, *where the supremum is taken over all measurable* $g\colon S \to \mathbf{R}$ *satisfying* $|g(x)| \le V(x)$ *for all* $x \in S$. *Formulas for* $\rho$ *and* $M$ *are given in Section* 2.1. *In particular,* $P^n g(x)$ *and* $\int g \, d\pi$ *are both well defined whenever*

$$\|g\|_V \equiv \sup\{|g(x)|/V(x)\colon x \in S\} < \infty.$$

The proof of Theorem 1.1 appears in Section 4. If we restrict to functions $g$ on the left-hand side of (1) that satisfy $|g(x)| \le 1$, we obtain the total variation norm $\|P_n(x, \cdot) - \pi\|_{\mathrm{TV}}$. So the inequality (1) is a strong version of the condition of *geometric ergodicity*, which says that for each $x \in S$ there exists $\gamma < 1$ such that

$$\gamma^{-n}\|P_n(x, \cdot) - \pi\|_{\mathrm{TV}} \to 0 \qquad \text{as } n \to \infty.$$

This concept was introduced in 1959 by Kendall [5] for countable state spaces. Important advances were made by Vere-Jones [22] in the countable setting, and by Nummelin and Tweedie [12] and Nummelin and Tuominen [11] for general state spaces. The condition in (1) is that of $V$-*uniform ergodicity*. Information about the theories of geometric ergodicity and $V$-uniform ergodicity is given in Chapters 15 and 16 of [8]. Results that relate the different notions of geometric ergodicity are also given in [13].

To date two basic methods have been used to obtain computable convergence rates. One method, introduced by Meyn and Tweedie [9], is based on *renewal theory*. In fact Theorem 1.1 is a restatement of Theorems 2.1–2.3 in [9], except that we give different formulas for $\rho$ and $M$. Our results in this paper use this method. The renewal theory method is easiest to describe when $C$ is an atom, that is, $P(x, A) = \nu(A)$ for all $x \in C$ and $A \in \mathcal{B}$. In this case, the Markov process $\{X_n : n \geq 0\}$ has a regeneration, or renewal, time whenever $X_n \in C$. Precise estimates are based on the regenerative decomposition, or first-entrance–last-exit decomposition; see the proof of Proposition 4.2. This method requires information about the regeneration time

$$\tau = \inf\{n > 0 : X_n \in C\}$$

which may be obtained using the drift condition (A2). It also requires information on the rate of convergence of the renewal sequence $u_n = P(X_n \in C | X_0 \in C)$ as $n \to \infty$. It is at this point that the aperiodicity condition (A3) is used. More generally, if $C$ is not an atom, then the renewal method may be applied to the split chain associated with the minorization condition (A1); see Section 4.2 for details of this construction.

The other main method, introduced by Rosenthal [18], is based on *coupling theory*, and relies on estimates of the coupling time $\hat{T} = \inf\{n > 0 : X_n = X_n'\}$ for some bivariate process $\{(X_n, X_n') : n \geq 0\}$ where each component is a copy of the original Markov chain. The minorization condition (A1) implies that the bivariate process can be constructed so that

$$P(X_{n+1} = X_{n+1}' | (X_n, X_n') \in C \times C) \geq \tilde{\beta}.$$

Therefore, coupling can be achieved with probability $\tilde{\beta}$ whenever $(X_n, X_n') \in C \times C$. It remains to estimate the hitting time $\inf\{n > 0 : (X_n, X_n') \in C \times C\}$. If the Markov chain is stochastically monotone and $C$ is a bottom or top set, then the univariate drift condition (A2) is sufficient. See the results in [7] and [21] for the case when $C$ is an atom, and in [16] for the general case. For stochastically monotone chains, the coupling method appears to be close to optimal. In the absence of stochastic monotonicity, a drift condition for the bivariate process is needed. This can often be achieved using the same function $V$ that appears in the (univariate) drift condition, but at the cost of enlarging the set $C$ and increasing the effective value of $\lambda$. Further information about these two methods and their relationship to our results appears in Section 7.

Our computations for $\rho$ and $M$ in Theorem 1.1 are valid for a very large class of Markov chains and, consequently, can be very far from sharp in particular cases. They can be improved dramatically in the setting of reversible Markov chains.

THEOREM 1.2. *Assume* (A1)–(A3) *and that the Markov chain is reversible (or symmetric) with respect to $\pi$, that is,*

$$\int_S Pf(x)g(x)\pi(dx) = \int_S f(x)Pg(x)\pi(dx)$$

*for all $f, g \in L^2(\pi)$. Then the assertions of Theorem* 1.1 *hold with the formulas for $\rho$ and $M$ in Section* 2.2.

Reversibility is an intrinsic feature of many MCMC algorithms, such as the Metropolis–Hastings algorithm and the random scan Gibbs sampler.

THEOREM 1.3. *In the setting of Theorem* 1.2 *assume also that the Markov chain is positive in the sense that*

$$\int_S Pf(x)f(x)\pi(dx) \geq 0$$

*for all $f \in L^2(\pi)$. Then the assertions of Theorem* 1.1 *hold with the formulas for $\rho$ and $M$ in Section* 2.3.

The proofs of Theorems 1.2 and 1.3 appear in Section 5, and some consequences for the spectral gap of $P$ in $L^2(\pi)$ appear in Section 6. For reversible positive Markov chains, our formulas for $\rho$ give the same values as the formulas given by Lund and Tweedie [7] (atomic case) and Roberts and Tweedie [16] (nonatomic case) under the assumption of stochastic monotonicity. The random scan Gibbs sampler is reversible and positive (see [6], Lemma 3). If $\{X_n : n \geq 0\}$ is one component in a two-component deterministic scan Gibbs sampler, then it is reversible and positive. Moreover, if a transition kernel $P$ is reversible with respect to $\pi$, then both the kernel $P^2$ for the two-skeleton chain and also the kernel $(I + P)/2$ for the binomial modification of the chain (see [20]) are reversible and positive. In particular, any discrete time skeleton of a continuous time reversible Markov process is positive.

In Section 8 we give numerical comparisons between our estimates and those obtained using [9] and the coupling method. The four Markov chains considered are "benchmark" examples used in earlier papers. Note that Theorem 1.1 outperforms the estimates given in [9]. For reversible chains, Theorem 1.2 is sometimes comparable with the coupling method, and sometimes noticeably better. For chains which are reversible and positive, Theorem 1.3 outperforms the coupling method.

In this paper our assumptions (A1)–(A3) all involve just the time 1 transition probabilities. In principle, our methods extend to a more general setting where one or more of the conditions involves $m$-step transitions for some $m > 1$. However, the calculations are much more cumbersome; we omit the details. Note that our method typically allows smaller $C$ than does the coupling method (see Section 7.2) and so there is less need to pass to minorization conditions involving time $m > 1$ (see the example in Section 8.4).

For the remainder of this introduction, we focus our attention on the formula for $\rho$. Define $\rho_V$ to be the infimum of all $\gamma$ for which an inequality of the form (1) holds true. Thus $\rho_V$ is the spectral radius of the operator $P -$ $1 \otimes \pi$ acting on the Banach space $(B_V, \| \cdot \|_V)$, say, of measurable functions $g : S \to \mathbf{R}$ such that $\|g\|_V < \infty$. We look for inequalities $\rho_V \leq \rho$, where $\rho$ is computable from the time 1 transition kernel.

At the heart of our calculations is an estimate on the rate of convergence of $\mathbf{P}^\nu(X_n \in C)$ to $\pi(C)$ as $n \to \infty$. More precisely, define

$$\rho_C = \limsup_{n \to \infty} |\mathbf{P}^\nu(X_n \in C) - \pi(C)|^{1/n}.$$

It is easy to verify [by taking $g(x) = \mathbb{1}_C(x)$ in (1), integrating with respect to $\nu$ and using $\int V \, d\nu < \infty$] that $\rho_C \leq \rho_V$. In the case that $C$ is an atom, we show (as a consequence of Propositions 4.1 and 4.2) that

$$(2) \qquad \rho_V \leq \max(\lambda, \rho_C).$$

Suppose instead that $C$ is not an atom, so that $\tilde{\beta} < 1$ in assumption (A1). We consider the associated split chain (see Section 4.2) and apply the atomic techniques to the split chain. In this case we show (as a consequence of Propositions 4.3 and 4.4) that

$$(3) \qquad \rho_V \leq \max(\lambda, (1 - \tilde{\beta})^{1/\alpha_1}, \rho_C),$$

where $\alpha_1 = 1 + (\log \frac{K - \tilde{\beta}}{1 - \tilde{\beta}})/(\log \lambda^{-1})$. We remark that $\max(\lambda, (1 - \tilde{\beta})^{1/\alpha_1}) = \beta_{\mathrm{RT}}^{-1}$, where $\beta_{\mathrm{RT}}$ is the estimate obtained by Roberts and Tweedie ([15], Theorem 2.3) for the radius of convergence of the generating function of the regeneration time for the split chain. Therefore, (3) may be rewritten $\rho_V \leq \max(\beta_{\mathrm{RT}}^{-1}, \rho_C)$.

It remains to get a good upper bound on $\rho_C$. We do this using renewal theory. Suppose first that $C$ is an atom and consider the renewal sequence $u_0 = 1$ and $u_n = \mathbf{P}(X_n \in C | X_0 \in C) = \mathbf{P}^\nu(X_{n-1} \in C)$ for $n \geq 1$. The $V$-uniform ergodicity implies that $\pi(C) = \lim_{n \to \infty} \mathbf{P}^\nu(X_{n-1} \in C) = \lim_{n \to \infty} u_n = u_\infty$, say. Thus $\rho_C^{-1}$ is the radius of convergence of the series $\sum_{n=1}^{\infty}(u_n - u_\infty)z^n$. The renewal sequence $u_n$, $n \geq 0$, is related to its corresponding increment sequence $b_n = \mathbf{P}^a(\tau = n)$, $n \geq 1$, by the renewal equation

$$u(z) = 1/(1 - b(z))$$

for $|z| < 1$, where $u(z) = \sum_{n=0}^{\infty} u_n z^n$ and $b(z) = \sum_{n=1}^{\infty} b_n z^n$. The drift condition (A2) implies that

$$\sum_{n=1}^{\infty} b_n \lambda^{-n} = \mathbf{E}^a(\lambda^{-\tau}) \leq \lambda^{-1} K$$

(see Proposition 4.1) and the aperiodicity condition (A3) implies that $b_1 = P(a, C) = \nu(C) \geq \beta$. In these circumstances a result of Kendall [5] shows that $\rho_C < 1$. In Section 3 we sharpen Kendall's result, using the lower bound on $b_1$ and the upper bound on $\sum_{n=1}^{\infty} b_n \lambda^{-n}$ to get an upper bound on $\rho_C$, depending only on $\lambda$, $K$ and $\beta$, which is strictly less than 1. In fact we give three different upper bounds on $\rho_C$. The first formula (in Theorem 3.2) is valid with no further restrictions on the Markov chain. The second formula (in Theorem 3.3) is valid for reversible Markov chains and the third formula (in Corollary 3.1) is valid for Markov chains which are reversible and positive.

The idea in the nonatomic case is similar. For the split chain the renewal sequence is given by $\bar{u}_n = \tilde{\beta}\mathbf{P}^\nu(X_{n-1} \in C)$ for $n \geq 1$, so that $\bar{u}_n \to \bar{u}_\infty$ has geometric convergence rate given by $\rho_C$. For the corresponding increment sequence $\bar{b}_n$, the estimate on $\sum_{n=1}^\infty \bar{b}_n r^n$ is more complicated, see (26) and (22), but the way in which results from Section 3 are applied is exactly the same.

**2. Formulas for $\rho$ and $M$.** Here we complete the statement of Theorems 1.1, 1.2 and 1.3 by giving formulas for the constants $\rho$ and $M$. We say that the set $C$ is an *atom* if $P(x, \cdot) = P(y, \cdot)$ for all $x, y \in C$. In this case we assume that $\tilde{\beta} = 1$ and $\nu = P(x, \cdot)$ for any $x \in C$. If $C$ is not an atom, so that $\tilde{\beta} < 1$, we define

$$\alpha_1 = 1 + \left(\log \frac{K - \tilde{\beta}}{1 - \tilde{\beta}}\right) \Big/ (\log \lambda^{-1})$$

and

$$\alpha_2 = 1 + \left(\log \frac{K}{\tilde{\beta}}\right) \Big/ (\log \lambda^{-1}).$$

In the special case when $\nu(C) = 1$, we can take $\alpha_2 = 1$. More generally, if we have the extra information that $\nu(C) + \int_{S \setminus C} V \, d\nu \leq \widetilde{K}$, we can take $\alpha_2 = 1 + (\log \widetilde{K})/(\log \lambda^{-1})$. Then define

$$R_0 = \min(\lambda^{-1}, (1 - \tilde{\beta})^{-1/\alpha_1})$$

and, for $1 < R \leq R_0$, define

$$L(R) = \frac{\tilde{\beta}R^{\alpha_2}}{1 - (1 - \tilde{\beta})R^{\alpha_1}}.$$

2.1. *Formulas for Theorem* 1.1. For $\beta > 0$, $R > 1$ and $L > 1$, define $R_1 = R_1(\beta, R, L)$ to be the unique solution $r \in (1, R)$ of the equation

$$(4) \qquad \frac{(r-1)}{r(\log R/r)^2} = \frac{e^2 \beta (R-1)}{8(L-1)}.$$

Since the left-hand side of (4) increases monotonically from 0 to $\infty$ as $r$ increases from 1 to $R$, the value $R_1$ is well defined and is easy to compute numerically. For $1 < r < R_1$, define

$$K_1(r, \beta, R, L) = \frac{2\beta + 2(\log N)(\log R/r)^{-1} - 8Ne^{-2}(r-1)r^{-1}(\log R/r)^{-2}}{(r-1)[\beta - 8Ne^{-2}(r-1)r^{-1}(\log R/r)^{-2}]},$$

where $N = (L-1)/(R-1)$.

*Atomic case.* We have $\rho = 1/R_1(\beta, \lambda^{-1}, \lambda^{-1}K)$ and, for $\rho < \gamma < 1$,

$$
\begin{aligned}
(5) \quad M = {} & \frac{\max(\lambda, K - \lambda/\gamma)}{\gamma - \lambda} + \frac{K(K - \lambda/\gamma)}{\gamma(\gamma - \lambda)} K_1(\gamma^{-1}, \beta, \lambda^{-1}, \lambda^{-1}K) \\
& + \frac{(K - \lambda/\gamma)\max(\lambda, K - \lambda)}{(\gamma - \lambda)(1 - \lambda)} + \frac{\lambda(K - 1)}{(\gamma - \lambda)(1 - \lambda)}.
\end{aligned}
$$

*Nonatomic case.* Let

$$
\widetilde{R} = \operatorname*{arg\,max}_{1 < R \le R_0} R_1(\beta, R, L(R)).
$$

Then $\rho = 1/R_1(\beta, \widetilde{R}, L(\widetilde{R}))$, and for $\rho < \gamma < 1$,

$$
\begin{aligned}
(6) \quad M = {} & \frac{\max(\lambda, K - \lambda/\gamma)}{\gamma - \lambda} + \frac{K[K\gamma - \lambda - \tilde{\beta}(\gamma - \lambda)]}{\gamma^2(\gamma - \lambda)[1 - (1 - \tilde{\beta})\gamma^{-\alpha_1}]} \\
& + \frac{\tilde{\beta}\gamma^{-\alpha_2 - 2}K(K\gamma - \lambda)}{(\gamma - \lambda)[1 - (1 - \tilde{\beta})\gamma^{-\alpha_1}]^2} K_1(\gamma^{-1}, \beta, \widetilde{R}, L(\widetilde{R})) \\
& + \frac{\gamma^{-\alpha_2 - 1}(K\gamma - \lambda)}{(\gamma - \lambda)[1 - (1 - \tilde{\beta})\gamma^{-\alpha_1}]^2} \\
& \quad \times \left( \frac{\tilde{\beta}\max(\lambda, K - \lambda)}{1 - \lambda} + \frac{(1 - \tilde{\beta})(\gamma^{-\alpha_1} - 1)}{\gamma^{-1} - 1} \right) \\
& + \frac{\gamma^{-\alpha_2}\lambda(K - 1)}{(1 - \lambda)(\gamma - \lambda)[1 - (1 - \tilde{\beta})\gamma^{-\alpha_1}]} \\
& + \frac{[K - \lambda - \tilde{\beta}(1 - \lambda)]}{(1 - \lambda)(1 - \gamma)} \left( (\gamma^{-\alpha_2} - 1) + \frac{(1 - \tilde{\beta})(\gamma^{-\alpha_1} - 1)}{\tilde{\beta}} \right).
\end{aligned}
$$

Notice that the result remains true with $\widetilde{R}$ replaced by any $R \in (1, R_0)$, but it does not give such a small $\rho$. We do not claim that $\widetilde{R}$ gives the smallest $K_1$.

2.2. *Formulas for Theorem* 1.2. Here we assume that the Markov chain is reversible.

*Atomic case.* Define

$$
R_2 = \begin{cases} \sup\{r < \lambda^{-1} : 1 + 2\beta r > r^{1 + (\log K)/(\log \lambda^{-1})}\}, & \text{if } K > \lambda + 2\beta, \\ \lambda^{-1}, & \text{if } K \le \lambda + 2\beta. \end{cases}
$$

Then $\rho = R_2^{-1}$ and, for $\rho < \gamma < 1$, replace $K_1(\gamma^{-1}, \beta, \lambda^{-1}, \lambda^{-1}K)$ by $K_2 = 1 + 1/(\gamma - \rho)$ in (5) for $M$ in Section 2.1. We remark that, using the convexity of $r^{1+(\log K)/(\log \lambda^{-1})}$, we can replace $\rho$ by the larger, but more easily computable, $\tilde{\rho}$ given by

$$\tilde{\rho} = \begin{cases} 1 - 2\beta(1-\lambda)/(K-\lambda), & \text{if } K > \lambda + 2\beta, \\ \lambda, & \text{if } K \leq \lambda + 2\beta. \end{cases}$$

*Nonatomic case.*  Define

$$R_2 = \begin{cases} \sup\{r < R_0 : 1 + 2\beta r > L(r)\}, & \text{if } L(R_0) > 1 + 2\beta R_0, \\ R_0, & \text{if } L(R_0) \leq 1 + 2\beta R_0. \end{cases}$$

Then $\rho = R_2^{-1}$ and, for $\rho < \gamma < 1$, replace $K_1(\gamma^{-1}, \beta, \widetilde{R}, L(\widetilde{R}))$ by $K_2 = 1 + \sqrt{\beta}/(\gamma - \rho)$ in (6) for $M$ given in Section 2.1.

2.3. *Formulas for Theorem* 1.3.  Here we assume that the Markov chain is reversible and positive.

*Atomic case.*  We have $\rho = \lambda$ and $M$ is calculated as in Section 2.2.

*Nonatomic case.*  We have $\rho = R_0^{-1}$ and $M$ is calculated as in Section 2.2.

**3. Kendall's theorem.**  The setting for this section is discrete renewal theory. Suppose that $V_1, V_2, \ldots$ are independent identically distributed random variables taking values in the set of positive integers and let $b_n = \mathbf{P}(V_1 = n)$ for $n \geq 1$. Define $T_0 = 0$ and $T_k = V_1 + \cdots + V_k$ for $k \geq 1$. Let $u_n = \mathbf{P}$ (there exists $k \geq 0$ such that $T_k = n$) for $n \geq 0$. Thus $u_n$ is the (undelayed) renewal sequence that corresponds to the increment sequence $b_n$. The following result is due to Kendall [5].

THEOREM 3.1.  *Assume that the sequence $\{b_n\}$ is aperiodic and that $\sum_{n=1}^{\infty} b_n R^n < \infty$ for some $R > 1$. Then $u_\infty = \lim_{n \to \infty} u_n$ exists and the series $\sum_{n=0}^{\infty} (u_n - u_\infty)z^n$ has radius of convergence greater than 1.*

In this section we obtain three different lower bounds on the radius of convergence of $\sum (u_n - u_\infty)z^n$.

3.1. *General case.*

THEOREM 3.2.  *Suppose that $\sum_{n=1}^{\infty} b_n R^n \leq L$ and $b_1 \geq \beta$ for some constants $R > 1$, $L < \infty$ and $\beta > 0$. Let $N = (L-1)/(R-1) \geq 1$. Let $R_1 = R_1(\beta, R, L)$ be the unique solution $r \in (1, R)$ of the equation*

$$\frac{(r-1)}{r(\log R/r)^2} = \frac{e^2\beta}{8N}.$$

*Then the series*

$$\sum_{n=1}^{\infty} (u_n - u_\infty) z^n$$

*has radius of convergence at least $R_1$. For any $r \in (1, R_1)$, define $K_1 = K_1(r, \beta, R, L)$ by*

$$K_1 = \frac{1}{r-1} \left( 1 + \frac{\beta + 2(\log N)(\log R/r)^{-1}}{\beta - 8Ne^{-2}(r-1)r^{-1}(\log R/r)^{-2}} \right).$$

*Then*

$$(7) \qquad \left| \sum_{n=0}^{\infty} (u_n - u_\infty) z^n \right| \leq K_1 \qquad \text{for all } |z| \leq r.$$

PROOF. Define the sequence $c_n = \sum_{k=n+1}^{\infty} b_k$ for $n \geq 0$ and define generating functions $b(z) = \sum_{n=1}^{\infty} b_n z^n$, $c(z) = \sum_{n=0}^{\infty} c_n z^n$ and $u(z) = \sum_{n=0}^{\infty} u_n z^n$ for $|z| < 1$. The renewal equation gives

$$(8) \qquad c(z) = \frac{1 - b(z)}{1 - z} = \frac{1}{(1-z)u(z)} = \frac{1}{1 - \sum_{n=1}^{\infty} (u_{n-1} - u_n) z^n}$$

for $|z| < 1$. Since the power series for $c(z)$ has nonnegative coefficients, for $|z| \leq R$ we have

$$|c(z)| \leq c(R) = \frac{b(R) - 1}{R - 1} \leq \frac{L - 1}{R - 1} = N$$

so that $c(z)$ is holomorphic on $|z| < R$. Now

$$\Re((1-z)c(z)) = \Re(1 - b(z))$$

$$= \sum_{n=1}^{\infty} b_n \Re(1 - z^n)$$

$$\geq \beta \Re(1 - z)$$

for $|z| \leq 1$. It follows that

$$|c(re^{i\theta})| \geq \beta \frac{\Re(1 - re^{i\theta})}{|1 - re^{i\theta}|} \geq \beta \left| \sin\left(\frac{\theta}{2}\right) \right|$$

for all $r \leq 1$. In particular, since $c(r) > 0$ for all $r \geq 0$, we see that $c(z) \neq 0$ whenever $|z| \leq 1$. For $1 \leq r < R$,

$$|c(re^{i\theta})| \geq \beta |\sin(\theta/2)| - |c(re^{i\theta}) - c(e^{i\theta})|$$

$$\geq \beta |\sin(\theta/2)| - (c(r) - c(1)).$$

Moreover, for $1 \leq r < R$,

$$\begin{aligned}
|c(re^{i\theta})| &\geq c(r) - |re^{i\theta} - r| \sup\{|c'(z)| : z \in [r, re^{i\theta}]\} \\
&\geq c(r) - |re^{i\theta} - r| c'(r) \\
&= c(r) - 2r|\sin(\theta/2)| c'(r).
\end{aligned}$$

Combining these two estimates we obtain

$$|c(re^{i\theta})| \geq \frac{\beta - A(r)}{\beta/c(r) + B(r)},$$

where $A(r) = 2rc'(r)[c(r) - c(1)]/c(r)$ and $B(r) = 2rc'(r)/c(r)$. Since the power series for $c$ has nonnegative coefficients, we may apply Hölder's inequality to obtain

$$c(s) \leq c(r) \left(\frac{s}{r}\right)^{(\log c(R)/c(r))/(\log R/r)}$$

for $0 < r < s < R$. Letting $s \searrow r$ gives

$$c'(r) \leq \frac{c(r)}{r} \frac{\log c(R)/c(r)}{\log R/r}$$

and, consequently,

$$c(r) - c(1) \leq \frac{(r-1)c(r)}{r} \frac{\log c(R)/c(r)}{\log R/r}$$

for $1 \leq r < R$. Thus we obtain the estimates

$$A(r) \leq \frac{2(r-1)}{r} c(r) \left[\log \frac{N}{c(r)}\right]^2 \left[\log \frac{R}{r}\right]^{-2}$$

and

$$B(r) \leq 2 \left[\log \frac{N}{c(r)}\right] \left[\log \frac{R}{r}\right]^{-1}.$$

Using the inequality

$$x \left[\log \frac{N}{x}\right]^2 \leq 4Ne^{-2}$$

for $0 < x < N$ in $A(r)$ and the inequality $c(r) \geq 1$ in $B(r)$ we get

$$A(r) \leq \frac{8Ne^{-2}(r-1)}{r} \left[\log \frac{R}{r}\right]^{-2}$$

and

$$B(r) \leq 2\log N \left[\log \frac{R}{r}\right]^{-1}.$$

Thus for $1 < r < R_1$ we have

$$(9) \qquad |c(re^{i\theta})| \geq \frac{\beta - 8Ne^{-2}(r-1)r^{-1}(\log R/r)^{-2}}{\beta + 2(\log N)(\log R/r)^{-1}} > 0.$$

Therefore $c(z) \neq 0$ for all $|z| < R_1$. Recalling (8), we see that $\sum_{n=1}^{\infty}(u_{n-1} - u_n)z^n$ is holomorphic on $|z| < R_1$ and, therefore, $r^n|u_{n-1} - u_n| \to 0$ as $n \to \infty$ for each $r < R_1$. It follows directly that $u_\infty = \lim_{n \to \infty} u_n$ exists and $r^n|u_n - u_\infty| \to 0$ as $n \to \infty$ for all $r < R_1$. Furthermore, using the fact $u_n - u_\infty = \sum_{m=n+1}^{\infty}(u_{m-1} - u_m)$, we get

$$\sum_{n=0}^{\infty}(u_n - u_\infty)z^n = \frac{1}{z-1}\left(\sum_{m=1}^{\infty}(u_{m-1} - u_m)z^m - (1 - u_\infty)\right)$$

whenever $1 < |z| < R_1$. Therefore, using (8) again, for $1 < r < R_1$ we have

$$\sup_{|z| \leq r}\left|\sum_{n=0}^{\infty}(u_n - u_\infty)z^n\right| = \sup_{|z|=r}\left|\sum_{n=0}^{\infty}(u_n - u_\infty)z^n\right| \leq \frac{1}{r-1}\left(1 + \sup_{|z|=r}\frac{1}{|c(z)|}\right)$$

and now (7) follows from (9). $\quad\square$

The estimates in Theorem 3.2 apply to a very general class of renewal sequences and as a result they are very far from the best possible in certain more restricted settings. We see in Theorem 3.3 and Corollary 3.1 that the estimates can be dramatically improved when we have extra information about the origin of the renewal sequence. Meanwhile, the following discussion shows that the estimate on the radius of convergence in Theorem 3.2 can be of the correct order of magnitude.

Suppose that $\beta$ and $L$ are fixed. Then as $R \searrow 1$ we have

$$(10) \qquad R_1 - 1 \sim \frac{e^2\beta}{8(L-1)}(R-1)^3.$$

The effect of the $(R-1)^3$ term is that, typically, $R_1$ is very much closer to 1 than $R$ is. This is a major contributing factor to the disappointing estimates obtained using Theorem 1.1 in the examples in Sections 8.1 and 8.2. However, in the absence of any further information beyond that given by the constants $\beta$, $R$ and $L$, the following calculations show that the term $(R-1)^3$ in (10) is optimal.

Consider the family of examples $b(z) = \beta z + (1-\beta)z^k$ for fixed $\beta$ and $k \to \infty$. For each $k$ there is a solution $z_k$ of the equation $\beta z + (1-\beta)z^k = 1$ near $e^{2\pi i/k}$. Calculating the asymptotic expansion for $z_k e^{-2\pi i/k}$ in powers of $1/k$ we obtain

$$z_k = e^{2\pi i/k}\left[1 - \left(\frac{2\pi\beta i}{1-\beta}\right)k^{-2} + \left(\frac{2\pi^2\beta}{(1-\beta)^2} + \frac{2\pi\beta^2 i}{(1-\beta)^2}\right)k^{-3} + O(k^{-4})\right]$$

and thus

$$|z_k| = 1 + \left(\frac{2\pi^2\beta}{(1-\beta)^2}\right)k^{-3} + O(k^{-4}).$$

For fixed $\beta$ and $L$ this example satisfies the conditions of Theorem 3.2 as long as $\beta R + (1-\beta)R^k = L$. As $k \to \infty$ we have $R - 1 \sim \log R \sim (1/k)\log\left(\frac{L-\beta}{1-\beta}\right)$ and thus

$$|z_k| - 1 \sim \left(\frac{2\pi^2\beta}{(1-\beta)^2}\right)\left[\log\left(\frac{L-\beta}{1-\beta}\right)\right]^{-3}(R-1)^3.$$

It is clear from the proof of Theorem 3.2 that any $r$ satisfying (7) must satisfy $r < |z_k|$. Thus the factor $(R-1)^3$ in (10) is optimal, although clearly the factor $e^2\beta/8(L-1)$ is not.

3.2. *Reversible case.* In this section we assume that the renewal sequence $u_n$ is generated by a Markov chain $\{X_n : n \geq 0\}$ which is reversible with respect to its invariant probability measure $\pi$. Thus

$$\pi(dx)P(x,dy) = \pi(dy)P(y,dx)$$

in the sense that the measures on $S \times S$ given by the left-hand and right-hand sides agree.

THEOREM 3.3. *Let $\{X_n : n \geq 0\}$ be a Markov chain which is reversible with respect to a probability measure $\pi$ and satisfies $P(x,dy) \geq \tilde{\beta}\mathbb{1}_C(x)\nu(dy)$ for some set $C$ and probability measure $\nu$. Let $\{u_n : n \geq 0\}$ be the renewal sequence given by $u_0 = 1$ and $u_n = \tilde{\beta}P^\nu(X_{n-1} \in C)$ for $n \geq 1$, and suppose that the corresponding increment sequence $\{b_n : n \geq 1\}$ satisfies $\sum_{n=1}^{\infty} b_n R^n \leq L$ and $b_1 \geq \beta$ for some constants $R > 1$, $L < \infty$ and $\beta > 0$. If $L > 1 + 2\beta R$, define $R_2 = R_2(\beta, R, L)$ to be the unique solution $r \in (1, R)$ of the equation*

$$1 + 2\beta r = r^{(\log L)/(\log R)}$$

*and let $R_2 = R$ otherwise. Then the series*

$$\sum_{n=0}^{\infty}(u_n - u_\infty)z^n$$

*has radius of convergence at least $R_2$. Moreover, if*

$$(11) \qquad \lim_{n \to \infty}\left|\int_C P^n\mathbb{1}_C(x)\pi(dx) - (\pi(C))^2\right|r^n < \infty \qquad \text{for all } r < R_2,$$

*then, for $1 < r < R_2$, we have*

$$(12) \qquad \sum_{n=1}^{\infty}|u_n - u_\infty|r^n \leq \frac{\sqrt{\beta}r}{1 - r/R_2}.$$

PROOF. Notice first that the discussion of split chains in Section 4.2 implies that $\{u_n : n \geq 0\}$ is indeed a renewal sequence. The reversibility implies that the transition operator $P$ for the original chain $\{X_n : n \geq 0\}$ acts as a self-adjoint contraction on the Hilbert space $L^2(\pi)$. We use $\langle \cdot, \cdot \rangle$ for the inner product in $L^2(\pi)$ and $\| \cdot \|$ for the corresponding norm. For any $A \subset S$ we have

$$\pi(A) = \int P(x, A) \pi(dx) \geq \tilde{\beta} \nu(A) \pi(C),$$

so that $\nu$ is absolutely continuous with respect to $\pi$ and has Radon–Nikodym derivative $d\nu/d\pi \leq 1/(\tilde{\beta}\pi(C))$. Throughout this proof we write $f = \mathbb{1}_C$ and $g = d\nu/d\pi$. Then $f, g \in L^2(\pi)$ with $\|f\|^2 = \pi(C)$ and $\|g\|^2 \leq 1/(\tilde{\beta}\pi(C))$. Now for $|z| < 1$,

$$(1 - z)u(z) = (1 - z) + \tilde{\beta}(1 - z) \sum_{n=1}^{\infty} \langle P^{n-1} f, g \rangle z^n$$

$$= (1 - z) + \tilde{\beta} z (1 - z) \langle (I - zP)^{-1} f, g \rangle.$$

Since $P$ is a self-adjoint contraction on $L^2(\pi)$, its spectrum is a subset of $[-1, 1]$ and we have a spectral resolution

$$P = \int \lambda \, dE(\lambda)$$

(see, e.g., [23], Section XI.6), where $E(1) = I$ and $\lim_{\lambda \nearrow -1} E(\lambda) = 0$. Write $F(\lambda) = \langle E(\lambda) f, g \rangle$. The function $F$ is of bounded variation and the corresponding signed measure $\mu_{f,g}$, say, is supported on $[-1, 1]$ and has total mass $|\mu_{f,g}|([-1, 1]) \leq \|f\| \cdot \|g\| \leq \tilde{\beta}^{-1/2}$. We obtain for $|z| < 1$,

$$(1 - z)u(z) = (1 - z) + \tilde{\beta} z (1 - z) \int_{[-1,1]} (1 - z\lambda)^{-1} \mu_{f,g}(d\lambda)$$

and so the function $(1 - z)u(z)$ has a holomorphic extension at least to $\{z \in \mathbf{C} : z^{-1} \notin [-1, 1]\} = \mathbf{C} \setminus ((-\infty, -1] \cup [1, \infty))$. The renewal equation gives

$$(1 - z)u(z) = \frac{1 - z}{1 - b(z)}$$

for $|z| < 1$, and the function $b$ is holomorphic in $B(0, R)$. It follows that the only solutions in $B(0, R)$ of the equation $b(z) = 1$ lie on one or the other of the intervals $(-R, -1]$ and $[1, R)$. Since $b'(1) > 0$, the zero of $b(z) - 1$ at $z = 1$ is a simple zero. For $1 < r \leq R$ we have $b(r) > b(1) = 1$. For $1 < r < R$ we also have $b(-r) \leq -2b_1 r + b(r)$. Using the estimate $b(r) \leq [b(R)]^{(\log r)/(\log R)} = r^{(\log L)/(\log R)}$, it follows that for $1 < r < R_2$ we have $b(-r) < 1$, where $R_2$ is given in the statement of the theorem. Thus $(1 - z)u(z)$ has a holomorphic extension to $B(0, R_2)$ and the first statement of the theorem follows as in the proof of Theorem 3.2.

Now we assume (11). Given $r < R_2$ we have

$$(13) \qquad |\langle P^n f, f \rangle - (\pi(C))^2| \leq Mr^{-n}$$

for some $M$ (depending on $r$). Recalling the spectral resolution, we have

$$\langle P^n f, f \rangle = \int_{[-1,1]} \lambda^n d\langle E(\lambda)f, f \rangle.$$

Letting $n \to \infty$ we get

$$\lim_{n\to\infty} \langle P^n f, f \rangle = \int_{\{1\}} d\langle E(\lambda)f, f \rangle$$

and so (13) may be rewritten as

$$(14) \qquad \left| \int_{[-1,1)} \lambda^n d\langle E(\lambda)f, f \rangle \right| \leq Mr^{-n}.$$

Now $\lambda \to \langle E(\lambda)f, f \rangle$ is an increasing function and hence corresponds to a positive measure $\mu_f$, say, on $[-1, 1]$. Letting $n \to \infty$ in (14) through the even integers, we see that $\mu_f([-1, -1/r)) = \mu_f((1/r, 1)) = 0$. This is true for all $r < R_2$ and so $\langle E(\lambda)f, f \rangle$ is constant on $[-1, -1/R_2]$ and on $(1/R_2, 1)$. It follows that $F(\lambda) = \langle E(\lambda)f, g \rangle$ is constant on these same intervals and so the support of $|\mu_{f,g}|$ is contained in $[-1/R_2, 1/R_2] \cup \{1\}$. Noting that

$$u_\infty = \tilde{\beta} \lim_{n\to\infty} \langle P^{n-1}f, g \rangle = \tilde{\beta} \lim_{n\to\infty} \int \lambda^{n-1} \mu_{f,g}(d\lambda) = \tilde{\beta} \mu_{f,g}(\{1\})$$

we get, for $n \geq 1$,

$$\begin{aligned}
|u_n - u_\infty| &= \tilde{\beta} \left| \int_{[-1/R_2, 1/R_2]} \lambda^{n-1} \mu_{f,g}(d\lambda) \right| \\
&\leq \tilde{\beta} \left( \frac{1}{R_2} \right)^{n-1} |\mu_{f,g}| \left( \left[ \frac{-1}{R_2}, \frac{1}{R_2} \right] \right) \\
&\leq \sqrt{\tilde{\beta}} \left( \frac{1}{R_2} \right)^{n-1}.
\end{aligned}$$

So for $r < R_2$, we get

$$\sum_{n=1}^{\infty} |u_n - u_\infty| r^n \leq \frac{\sqrt{\tilde{\beta}} r}{1 - r/R_2}$$

as required.   $\square$

REMARK 3.1. The estimate (12) is true without the extra assumption (11) if $P$ is a compact operator on $L^2(\pi)$. The first assertion in Theorem 3.3 implies that

$$r^n \int_{[-1,1)} \lambda^n \mu_{f,g}(d\lambda) \to 0 \qquad \text{as } n \to \infty$$

for all $r < R_2$ and the compactness implies that the restriction of $\mu_{f,g}$ to $[-1, 1] \setminus [-1/R_2, 1/R_2]$ is a finite sum of atoms. It then follows directly that the support of $|\mu_{f,g}|$ is contained in $[-1/R_2, 1/R_2] \cup \{1\}$.

COROLLARY 3.1. *In the setting of Theorem* 3.3, *assume also that*

$$\int Pf(x)f(x)\pi(dx) \geq 0 \qquad \text{for all } f \in L^2(\pi).$$

*Then in the assertions of Theorem* 3.3 *we can take* $R_2 = R$.

PROOF. The additional assumption implies that the spectrum of $P$ is contained in $[0, 1]$. Arguing as in the proof of Theorem 3.3, we obtain, for $|z| < 1$,

$$(1-z)u(z) = (1-z) + \tilde{\beta}z(1-z)\int_{[0,1]}(1-z\lambda)^{-1}\mu_{f,g}(d\lambda)$$

and so the function $(1-z)u(z)$ has a holomorphic extension at least to $\{z \in \mathbf{C} : z^{-1} \notin [0, 1]\} = \mathbf{C} \setminus [1, \infty)$. It follows that the equation $b(z) = 1$ cannot have a solution in $(-R, -1]$ and so $(1-z)u(z)$ is holomorphic on $B(0, R)$. The remainder of the proof goes as in Theorem 3.3. $\square$

The following lemma enables us to apply Corollary 3.1 and Theorem 1.3 to a large class of Metropolis–Hastings chains, including the example in Section 8.2.

LEMMA 3.1. *The Metropolis–Hastings chain generated by a candidate transition density* $q(x, y)$ *of the form*

$$q(x, y) = \int r(z, x)r(z, y)\, dz$$

*is reversible and positive.*

PROOF. Since a Metropolis–Hastings chain is automatically reversible, it suffices to check positivity. For notational convenience, we identify the measure $\pi$ with its density $\pi(x)$ with respect to the reference measure $dx$. Notice first that for any $g \in L^2(\pi)$ we have

$$\iint g(x)g(y)\min(\pi(x), \pi(y))\, dx\, dy$$

$$= \iint g(x)g(y)\left(\int_0^\infty \mathbb{1}_{[0,\pi(x)]}(t)\mathbb{1}_{[0,\pi(y)]}(t)\right)dt\right)dx\, dy$$

$$(15) \qquad = \int_0^\infty \left(\iint g(x)\mathbb{1}_{[0,\pi(x)]}(t)g(y)\mathbb{1}_{[0,\pi(y)]}(t)\, dx\, dy\right)dt$$

$$= \int_0^\infty \left( \int g(x) \mathbb{1}_{[0,\pi(x)]}(t)\, dx \right)^2 dt$$

$$\geq 0.$$

The assumption on $q$ implies that $q(x, y) = q(y, x)$, and so the kernel $P$ for the Metropolis–Hastings chain is given by

$$Pf(x) = \int f(y) \min(\pi(y)/\pi(x), 1) q(x, y)\, dy + \alpha(x)f(x)$$

for some $\alpha(x) \geq 0$. Then, for $f \in L^2(\pi)$, we have

$$\int Pf(x)f(x)\pi(x)\, dx = \iint f(x)f(y) \min(\pi(x), \pi(y)) q(x, y)\, dx\, dy$$
$$+ \int \alpha(x)f(x)^2 \pi(x)\, dx.$$

Clearly the second term on the right-hand side is nonnegative, and the first term on the right-hand side is

$$\iint f(x)f(y) \min(\pi(x), \pi(y)) \left( \int r(z, x)r(z, y)\, dz \right) dx\, dy$$
$$= \int \left( \iiint f(x)r(z, x)f(y)r(z, y) \min(\pi(x), \pi(y))\, dx\, dy \right) dz$$
$$\geq 0,$$

where we use (15) with $g(x) = f(x)r(z, x)$ and then integrate with respect to $z$. □

REMARK 3.2. The condition on $q$ is satisfied if $r$ is a symmetric Markov kernel and $q$ corresponds to two steps of $r$.

**4. Proof of Theorem 1.1.** In this section we describe the methods used to obtain the formulas in Section 2.1 for $\rho$ and $M$. From the results of Meyn and Tweedie [8, 9] we know that $\{X_n : n \geq 0\}$ is $V$-uniformly ergodic, with invariant probability measure $\pi$, say. We concentrate on the calculation of $\rho$ and $M$. We do not make any assumption of reversibility in this section. At the appropriate point in the argument we appeal to Theorem 3.2. Proofs of Propositions 4.1–4.4 appear in the Appendix.

4.1. *Atomic case.* Suppose that $C$ is an atom for the Markov chain. Then in the minorization condition (A1) we can take $\tilde{\beta} = 1$ and $\nu = P(a, \cdot)$ for some fixed point $a \in C$. Let $\tau$ be the stopping time

$$\tau = \inf\{n \geq 1 : X_n \in C\}$$

and define $u_n = \mathbf{P}^a(X_n \in C)$ for $n \geq 0$. Then $u_n$ is the renewal sequence that corresponds to the increment sequence $b_n = \mathbf{P}^a(\tau = n)$ for $n \geq 1$. Define functions $G(r, x)$ and $H(r, x)$ by

$$G(r, x) = \mathbf{E}^x(r^\tau),$$

$$H(r, x) = \mathbf{E}^x\left(\sum_{n=1}^{\tau} r^n V(X_n)\right)$$

for all $x \in S$ and all $r > 0$ for which the right-hand sides are defined. Most of the following result is well known (see, e.g., [7], Lemma 2.2 and Theorem 3.1). The estimate in (iv) appears to be new, and helps to reduce our estimate for $M$.

PROPOSITION 4.1. *Assume only the drift condition* (A2).

(i) *For all* $x \in S$, $\mathbf{P}^x(\tau < \infty) = 1$.

(ii) *For* $1 \leq r \leq \lambda^{-1}$,

$$G(r, x) \leq \begin{cases} V(x), & \text{if } x \notin C, \\ rK, & \text{if } x \in C. \end{cases}$$

(iii) *For* $0 < r < \lambda^{-1}$,

$$H(r, x) \leq \begin{cases} \dfrac{r\lambda V(x)}{1 - r\lambda}, & \text{if } x \notin C, \\ \dfrac{r(K - r\lambda)}{1 - r\lambda}, & \text{if } x \in C. \end{cases}$$

(iv) *For* $1 < r < \lambda^{-1}$ *and* $x \in C$,

$$\frac{H(r, x) - rH(1, x)}{r - 1} \leq \frac{\lambda r(K - 1)}{(1 - \lambda)(1 - r\lambda)}.$$

The following result is a minor variation of results in [8].

PROPOSITION 4.2. *Assume only that the Markov chain is geometrically ergodic with (unique) invariant probability measure* $\pi$, *that* $C$ *is an atom and that* $V$ *is a nonnegative function. Suppose* $g : S \to \mathbf{R}$ *satisfies* $\|g\|_V \leq 1$. *Then*

$$\sup_{|z| \leq r} \left| \sum_{n=1}^{\infty} \left( P^n g(x) - \int g \, d\pi \right) z^n \right|$$

$$\leq H(r, x) + G(r, x)H(r, a) \sup_{|z| \leq r} \left| \sum_{n=0}^{\infty} (u_n - u_\infty) z^n \right|$$

$$+ H(r, a) \frac{G(r, x) - 1}{r - 1} + \frac{H(r, a) - rH(1, a)}{r - 1}$$

*for all* $r > 1$ *for which the right-hand side is finite.*

It is an immediate consequence of Propositions 4.1 and 4.2 that $\rho_V \le \max(\lambda, \rho_C)$ when $C$ is an atom.

PROOF OF ESTIMATES FOR THE ATOMIC CASE.   We apply Theorem 3.2 to the sequence $u_n$. For the increment sequence $b_n = P^a(\tau = n)$ we have $\sum_{n=1}^{\infty} b_n \lambda^{-n} = \mathbf{E}^a(\lambda^{-\tau}) = G(\lambda^{-1}, a) \le \lambda^{-1} K$. Moreover the aperiodicity condition (A3) gives $b_1 = P(a, C) \ge \beta$. For $1 < r < R_1(\beta, \lambda^{-1}, \lambda^{-1} K)$ and $K_1 = K_1(r, \beta, \lambda^{-1}, \lambda^{-1} K)$, Theorem 3.2 gives

$$\sup_{|z| \le r} \left| \sum_{n=0}^{\infty} (u_n - u_\infty) z^n \right| \le K_1.$$

By substituting this and the estimates from Proposition 4.1 into Proposition 4.2 together with the inequality

$$\frac{G(r, x) - 1}{r - 1} \le \frac{G(\lambda^{-1}, x) - 1}{\lambda^{-1} - 1} \le \frac{\max(\lambda, K - \lambda)}{1 - \lambda} V(x),$$

we get

$$\sup_{|z| \le r} \left| \sum_{n=1}^{\infty} \left( P^n g(x) - \int g \, d\pi \right) z^n \right| \le M V(x)$$

and so

$$\left| P^n g(x) - \int g \, d\pi \right| \le M V(x) r^{-n},$$

where

$$(16) \quad \begin{aligned} M &= \frac{r \max(\lambda, K - r\lambda)}{1 - r\lambda} + \frac{r^2 K(K - r\lambda)}{(1 - r\lambda)} K_1(r, \beta, \lambda^{-1}, \lambda^{-1} K) \\ &\quad + \frac{r(K - r\lambda) \max(\lambda, K - \lambda)}{(1 - r\lambda)(1 - \lambda)} + \frac{\lambda r(K - 1)}{(1 - r\lambda)(1 - \lambda)}. \end{aligned}$$

Therefore, we can take $\rho = 1/R_1(\beta, \lambda^{-1}, \lambda^{-1} K)$ and the formula for $M$ is obtained by putting $r = 1/\gamma$ in (16).   □

4.2. *Nonatomic case.*   If $C$ is not an atom, then in the minorization condition (A1) we must have $\tilde{\beta} < 1$. Following Nummelin ([10], Section 4.4), we consider the split chain $\{(X_n, Y_n) : n \ge 0\}$ with state space $S \times \{0, 1\}$ and transition probabilities given by

$$P\{Y_n = 1 | \mathcal{F}_n^X \vee \mathcal{F}_{n-1}^Y\} = \tilde{\beta} \mathbb{1}_C(X_n),$$

$$P\{X_{n+1} \in A | \mathcal{F}_n^X \vee \mathcal{F}_n^Y\} = \begin{cases} \nu(A), & \text{if } Y_n = 1, \\ \dfrac{P(X_n, A) - \tilde{\beta} \mathbb{1}_C(X_n) \nu(A)}{1 - \tilde{\beta} \mathbb{1}_C(X_n)}, & \text{if } Y_n = 0. \end{cases}$$

Here $\mathcal{F}_n^X = \sigma\{X_r : 0 \leq r \leq n\}$ and $\mathcal{F}_n^Y = \sigma\{Y_r : 0 \leq r \leq n\}$. Thus the split chain evolves as follows. Given $X_n$, choose $Y_n$ so that $\mathbf{P}(Y_n = 1) = \tilde{\beta}\mathbb{1}_C(X_n)$. If $Y_n = 1$ then $X_{n+1}$ has distribution $\nu$, whereas if $Y_n = 0$ then $X_{n+1}$ has distribution $(P(X_n, \cdot) - \tilde{\beta}\mathbb{1}_C(X_n)\nu)/(1 - \tilde{\beta}\mathbb{1}_C(X_n))$. The split chain $\{(X_n, Y_n) : n \geq 0\}$ is designed so that it has an atom $S \times \{1\}$ and so that its first component $\{X_n : n \geq 0\}$ is a copy of the original Markov chain.

We apply the ideas of Section 4.1 to the split chain $(X_n, Y_n)$ with atom $S \times \{1\}$ and stopping time

(17) $$T = \min\{n \geq 1 : Y_n = 1\}.$$

Let $\overline{\mathbf{P}}^{x,i}$ and $\overline{\mathbf{E}}^{x,i}$ denote probability and expectation for the split chain started with $X_0 = x$ and $Y_0 = i$. To emphasize the similarities with the calculations in the previous section, we fix a point $a \in C$, and write $\overline{\mathbf{P}}^{x,1} = \overline{\mathbf{P}}^{a,1}$ and $\overline{\mathbf{E}}^{x,1} = \overline{\mathbf{E}}^{a,1}$. Define the renewal sequence $\bar{u}_n = \overline{\mathbf{P}}^{a,1}(Y_n = 1)$ for $n \geq 0$ and the corresponding increment sequence $\bar{b}_n = \overline{\mathbf{P}}^{a,1}(T = n)$ for $n \geq 1$. Notice that $\bar{u}_n = \tilde{\beta}\overline{\mathbf{P}}^{a,1}(X_n \in C) = \tilde{\beta}\mathbf{P}^\nu(X_{n-1} \in C)$ for $n \geq 1$, so that $\rho_C$ controls the rate of convergence of $\bar{u}_n \to \bar{u}_\infty$ in the nonatomic case also. Following the methods used in the atomic case, we define

$$\overline{G}(r, x, i) = \overline{\mathbf{E}}^{x,i}(r^T),$$

$$\overline{H}(r, x, i) = \overline{\mathbf{E}}^{x,i}\left(\sum_{n=1}^{T} r^n V(X_n)\right)$$

for all $x \in S$, $i = 0, 1$ and all $r > 0$ for which the right-hand sides are defined. If we define

$$\overline{\mathbf{E}}^x = [1 - \tilde{\beta}\mathbb{1}_C(x)]\overline{\mathbf{E}}^{x,0} + \tilde{\beta}\mathbb{1}_C(x)\overline{\mathbf{E}}^{x,1},$$

then $\overline{\mathbf{E}}^x$ agrees with $\mathbf{E}^x$ on $\mathcal{F}^X = \sigma\{X_n : n \geq 0\}$. Define

$$\overline{G}(r, x) = \overline{\mathbf{E}}^x(r^T),$$

$$\overline{H}(r, x) = \overline{\mathbf{E}}^x\left(\sum_{n=1}^{T} r^n V(X_n)\right).$$

Applying the techniques used in Proposition 4.2 to the split chain, we obtain the following result.

PROPOSITION 4.3. *Assume only that the original Markov chain is geometrically ergodic with (unique) invariant probability measure $\pi$ and that $V$ is a nonnegative function. Suppose $g : S \to \mathbf{R}$ satisfies $\|g\|_V \leq 1$. Then*

$$\sup_{|z| \leq r}\left|\sum_{n=1}^{\infty}\left(P^n g(x) - \int g\, d\pi\right)z^n\right|$$

$$\leq \overline{H}(r,x) + \overline{G}(r,x)\overline{H}(r,a,1) \sup_{|z|\leq r}\left|\sum_{n=0}^{\infty}(\bar{u}_n - \bar{u}_\infty)z^n\right|$$

$$+ \overline{H}(r,a,1)\frac{\overline{G}(r,x)-1}{r-1} + \frac{\overline{H}(r,a,1)-r\overline{H}(1,a,1)}{r-1}$$

for all $r > 1$ for which the right-hand side is finite.

We need to extend the estimates on $G(r,x)$ and $H(r,x)$ from Section 4.1 to estimates on the corresponding functions $\overline{G}(r,x,i)$ and $\overline{H}(r,x,i)$ defined in terms of the split chain and the stopping time $T$. Define

$$\widetilde{G}(r) = \sup\{\overline{\mathbf{E}}^{x,0}(r^\tau) : x \in C\}.$$

Notice that the initial condition $(x,0)$ for $x \in C$ represents a failed opportunity for the split chain to renew. Thus $\widetilde{G}(r)$ represents the extra contribution to $\overline{G}(r,x,i)$ and $\overline{H}(r,x,i)$ which occurs every time the split chain has $X_n \in C$ but fails to have $Y_n = 1$. Given $X_n \in C$, this failure occurs with probability $(1-\tilde{\beta})$. Thus to get finite estimates for $\overline{G}(r,x,i)$ and $\overline{H}(r,x,i)$, we insist on the condition $(1-\tilde{\beta})\widetilde{G}(r) < 1$. This idea is formalized in Lemmas A.1 and A.2 in the Appendix. For our purposes here the important estimates are given in the following result. The estimate (19) and an estimate closely related to (21) appear in [15], where they denote $R_0 = \beta_{\mathrm{RT}}$.

PROPOSITION 4.4.   *Assume conditions* (A1) *and* (A2) *with* $\tilde{\beta} < 1$. *Define*

$$(18) \qquad \alpha_1 = 1 + \left(\log\frac{K-\tilde{\beta}}{1-\tilde{\beta}}\right)\Big/(\log\lambda^{-1}).$$

*Then, for* $1 \leq r \leq \lambda^{-1}$,

$$(19) \qquad\qquad\qquad \widetilde{G}(r) \leq r^{\alpha_1}.$$

*Furthermore, define*

$$(20) \qquad\qquad \alpha_2 = 1 + \left(\log\frac{K}{\tilde{\beta}}\right)\Big/(\log\lambda^{-1})$$

*and* $R_0 = \min(\lambda^{-1}, (1-\tilde{\beta})^{-1/\alpha_1})$. *Then*

$$(21) \qquad\qquad \overline{G}(r,x) \leq \frac{\tilde{\beta}G(r,x)}{1-(1-\tilde{\beta})r^{\alpha_1}},$$

$$(22) \qquad\qquad \overline{G}(r,a,1) \leq \frac{\tilde{\beta}r^{\alpha_2}}{1-(1-\tilde{\beta})r^{\alpha_1}} \equiv L(r),$$

$$(23) \qquad\qquad \overline{H}(r,x) \leq H(r,x) + \frac{r[K-r\lambda-\tilde{\beta}(1-r\lambda)]}{(1-r\lambda)[1-(1-\tilde{\beta})r^{\alpha_1}]}G(r,x),$$

$$(24) \qquad \overline{H}(r,a,1) \leq \frac{r^{\alpha_2+1}(K-r\lambda)}{(1-r\lambda)[1-(1-\tilde{\beta})r^{\alpha_1}]},$$

$$\frac{\overline{H}(r,a,1) - r\overline{H}(1,a,1)}{r-1}$$

$$(25) \qquad \leq \frac{r^{\alpha_2+1}\lambda(K-1)}{(1-\lambda)(1-r\lambda)[1-(1-\tilde{\beta})r^{\alpha_1}]}$$

$$+ \frac{r[K-\lambda-\tilde{\beta}(1-\lambda)]}{(1-\lambda)[1-(1-\tilde{\beta})r^{\alpha_1}]}\left(\frac{r^{\alpha_2}-1}{r-1} + \frac{(1-\tilde{\beta})(r^{\alpha_1}-1)}{\tilde{\beta}(r-1)}\right)$$

whenever $1 < r < R_0$.

REMARK 4.1. If $\nu(C) = 1$, then $G(r,a,1) = r$ and so we can take $\alpha_2 = 1$ in Proposition 4.4. More generally if we know that $\nu(C) + \int_{S \setminus C} V \, d\nu \leq \widetilde{K}$, then we can take $\alpha_2 = 1 + (\log \widetilde{K})/(\log \lambda^{-1})$.

It is an immediate consequence of Propositions 4.3 and 4.4 that

$$\rho_V \leq \max(\lambda, \rho_C, (1-\tilde{\beta})^{1/\alpha_1})$$

when $C$ is not an atom.

PROOF OF ESTIMATES FOR THE NONATOMIC CASE. We apply Theorem 3.2 to the sequence $\bar{u}_n$. For the increment sequence $\bar{b}_n = \overline{\mathbf{P}}^{a,1}(T = n)$ we have

$$(26) \qquad \sum_{n=1}^{\infty} \bar{b}_n R^n = \overline{\mathbf{E}}^{a,1}(R^T) = \overline{G}(R,a,1) \leq L(R)$$

for $1 < R < R_0$, where the constant $R_0$ and the function $L(R)$ are defined in Proposition 4.4. The aperiodicity condition (A3) implies $\bar{b}_1 = \tilde{\beta}P(a,C) \geq \beta$. For the moment fix a value of $R$ in the range $1 < R < R_0$. By Theorem 3.2, for $1 < r < R_1(\beta, R, L(R))$, we have

$$\sup_{|z| \leq r}\left|\sum_{n=0}^{\infty}(\bar{u}_n - \bar{u}_\infty)z^n\right| \leq K_1(r, \beta, R, L(R)).$$

Notice that (21) implies

$$\frac{\overline{G}(r,x)-1}{r-1} \leq \frac{1}{1-(1-\tilde{\beta})r^{\alpha_1}}\left[\tilde{\beta}\left(\frac{G(r,x)-1}{r-1}\right) + (1-\tilde{\beta})\left(\frac{r^{\alpha_1}-1}{r-1}\right)\right].$$

Then using the estimates from Propositions 4.1 and 4.4 in Proposition 4.3 we get, for $1 < r < R_1(\beta, R, L(R))$,

$$\sup_{|z| \leq r} \left| \sum_{n=1}^{\infty} \left( P^n g(x) - \int g \, d\pi \right) z^n \right| \leq MV(x),$$

where

$$
\begin{aligned}
M = {} & \frac{r \max(\lambda, K - r\lambda)}{1 - r\lambda} + \frac{r^2 K[K - r\lambda - \tilde{\beta}(1 - r\lambda)]}{(1 - r\lambda)[1 - (1 - \tilde{\beta})r^{\alpha_1}]} \\
& + \frac{\tilde{\beta} r^{\alpha_2 + 2} K(K - r\lambda)}{(1 - r\lambda)[1 - (1 - \tilde{\beta})r^{\alpha_1}]^2} K_1(r, \beta, R, L(R)) \\
& + \frac{r^{\alpha_2 + 1}(K - r\lambda)}{(1 - r\lambda)[1 - (1 - \tilde{\beta})r^{\alpha_1}]^2} \\
& \quad \times \left( \frac{\tilde{\beta} \max(\lambda, K - \lambda)}{1 - \lambda} + \frac{(1 - \tilde{\beta})(r^{\alpha_1} - 1)}{r - 1} \right) \\
& + \frac{r^{\alpha_2 + 1} \lambda(K - 1)}{(1 - \lambda)(1 - r\lambda)[1 - (1 - \tilde{\beta})r^{\alpha_1}]} \\
& + \frac{r[K - \lambda - \tilde{\beta}(1 - \lambda)]}{(1 - \lambda)[1 - (1 - \tilde{\beta})r^{\alpha_1}]} \left( \frac{r^{\alpha_2} - 1}{r - 1} + \frac{(1 - \tilde{\beta})(r^{\alpha_1} - 1)}{\tilde{\beta}(r - 1)} \right).
\end{aligned}
$$

(27)

To obtain the smallest possible $\rho$, we choose $\widetilde{R} \in (1, R_0)$ so as to maximize $R_1(\beta, R, L(R))$. Then we take $\rho = 1/R_1(\beta, \widetilde{R}, L(\widetilde{R}))$ and substitute $r = \gamma^{-1}$ in formula (27) for $M$ and we are done. $\square$

## 5. Proof of Theorems 1.2 and 1.3

In this section we assume that the Markov chain $\{X_n : n \geq 0\}$ is reversible with respect to its invariant probability measure $\pi$. We first obtain the estimates of Section 2.2.

*Atomic case.* The proof in Section 4.1 goes through up to the point where we apply Theorem 3.2. Since the Markov chain is reversible, we can replace $R_1(\beta, \lambda^{-1}, \lambda^{-1} K)$ of Theorem 3.2 by $R_2 = R_2(\beta, \lambda^{-1}, \lambda^{-1} K)$ of Theorem 3.3. Then by the first part of Theorem 3.3, for $1 < r < R_2$ we, have

$$\sup_{|z| \leq r} \left| \sum_{n=0}^{\infty} (u_n - u_\infty) z^n \right| \leq K_2$$

for some $K_2 < \infty$. At this point we do not have an estimate for $K_2$. Continuing as in Section 4.1, we obtain, for $1 < r < R_2$,

$$\sup_{|z| \leq r} \left| \sum_{n=1}^{\infty} \left( P^n g(x) - \int g \, d\pi \right) z^n \right| \leq MV(x)$$

(28)

for some constant $M$. At this point we do not have an estimate for $M$. However, now in (28) we can take $g = \mathbb{1}_C$ and integrate the $x$ variable with respect to $\pi$ over $C$ to obtain the estimate (11). We can now apply the second part of Theorem 3.3 to obtain $K_2 = 1 + r/(1 - r/R_2)$. The rest of the proof goes as in Section 4.1. We have $\rho = 1/R_2(\beta, \lambda^{-1}, \lambda^{-1}K)$ and in (16) for $M$ we replace $K_1$ by $K_2$.

*Nonatomic case.* We have the estimate $\sum_{n=1}^{\infty} \bar{b}_n R^n = \overline{G}(R, a, 1) \le L(R)$ valid for all $1 \le R < R_0$, and we can choose the $R$ for which we apply Theorem 3.3. If $1 + 2\beta R_0 \ge L(R_0)$, then $L(R_0) < \infty$ and $\sum_{n=1}^{\infty} \bar{b}_n R_0^n = \overline{G}(R_0, a, 1) \le L(R_0)$. We can apply Theorem 3.3 with $R = R_0$ and obtain $R_2 = R_0$. This case can occur only when $R_0 = \lambda^{-1} < (1 - \tilde{\beta})^{-1/\alpha_1}$. Otherwise we take $\widetilde{R}$ to be the unique solution in the interval $(1, R_0)$ of the equation $1 + 2\beta R = L(R)$ and apply Theorem 3.3 with $R = \widetilde{R}$ to obtain $R_2 = \widetilde{R}$. Then by the first part of Theorem 3.3, for $1 < r < R_2$, we have

$$\sup_{|z| \le r} \left| \sum_{n=0}^{\infty} (\bar{u}_n - \bar{u}_\infty) z^n \right| \le K_2$$

for some $K_2 < \infty$. Initially we do not have an estimate for $K_2$, but the same method as above allows us to use the second part of Theorem 3.3 and assert that $K_2 = 1 + \sqrt{\beta} r/(1 - r/R_2)$. The rest of the proof goes as in Section 4.2. We have $\rho = 1/R_2$, where $R_2 = \sup\{r < R_0 : 1 + 2\beta r \ge L(r)\}$, and in (27) for $M$ we replace $K_1$ by $K_2$.

The estimates of Section 2.3 are obtained in a similar manner, using Corollary 3.1 in place of Theorem 3.3.

## 6. $L^2$-geometric ergodicity for reversible chains.

When the Markov chain is reversible with respect to the probability measure $\pi$, the Markov operator $P$ acts as a self-adjoint operator on $L^2(\pi)$. The equivalence of (pointwise) geometric ergodicity and the existence of a spectral gap for $P$ acting on $L^2(\pi)$ was proved in [13]. Also see [17] for the equivalence of $L^2$- and $L^1$-geometric ergodicity for reversible Markov chains.

THEOREM 6.1. *Assume that the Markov chain $\{X_n : n \ge 0\}$ is $V$-uniformly ergodic with invariant probability $\pi$ (so that $\int V \, d\pi < \infty$) and let $\rho_V$ be the spectral radius of $P - 1 \otimes \pi$ on $B_V$. Suppose also that $\{X_n : n \ge 0\}$ is reversible with respect to $\pi$. Then, for all $f \in L^2(\pi)$, we have*

$$\left\| P^n f - \int f \, d\pi \right\|_{L^2} \le (\rho_V)^n \left\| f - \int f \, d\pi \right\|_{L^2}.$$

*In particular, the spectral radius of $P - 1 \otimes \pi$ on $L^2(\pi)$ is at most $\rho_V$.*

PROOF. For ease of notation write $\int f \, d\pi = \bar{f}$. Suppose first that $f$ is a bounded function, so that $\|f\|_V < \infty$ and $\int |f(x)|V(x) \, d\pi(x) < \infty$. For any $\gamma > \rho_V$ there is $M < \infty$ so that

$$|P^n f(x) - \bar{f}| \le M\|f\|_V V(x)\gamma^n.$$

Multiplying by $f(x)$ and integrating with respect to $\pi$ we get

$$|\langle P^n f, f \rangle - \bar{f}^2| \le M\|f\|_V \left( \int |f(x)|V(x) \, d\pi(x) \right)\gamma^n.$$

Arguing as in the proof of Theorem 3.3 we see that for any $g \in L^2(\pi)$ the function $\lambda \mapsto \langle E(\lambda)f, g \rangle$ is constant on $[-1, -\rho_V)$ and on $(\rho_V, 1)$. The corresponding signed measure $\mu_{f,g}$ has $|\mu_{f,g}|([-1,1]) \le \|f\|_{L^2}\|g\|_{L^2}$. Therefore

$$|\langle P^n f - \bar{f}, g \rangle| = \left| \int_{[-\rho_V, \rho_V]} \lambda^n \, d\mu_{f,g}(\lambda) \right| \le \rho_V^n \|f\|_{L^2}\|g\|_{L^2}.$$

This is true for all $g \in L^2(\pi)$ so we obtain $\|P^n f - \bar{f}\|_{L^2} \le \rho_V^n \|f\|_{L^2}$. Replacing $f$ by $f - \bar{f}$ we obtain $\|P^n f - \bar{f}\|_{L^2} \le \rho_V^n \|f - \bar{f}\|_{L^2}$. Finally for arbitrary $f \in L^2(\pi)$ there exist bounded $f_k$ so that $\|f - f_k\|_{L^2} \to 0$. Then, for each $n \ge 0$,

$$\|P^n f - \bar{f}\|_{L^2} = \lim_{k \to \infty} \|P^n f_k - \bar{f}_k\|_{L^2}$$

$$\le \rho_V^n \lim_{k \to \infty} \|f_k - \bar{f}_k\|_{L^2} = \rho_V^n \|f - \bar{f}\|_{L^2}$$

and we are done. □

COROLLARY 6.1. *Assume that the Markov chain $\{X_n : n \ge 0\}$ satisfies* (A1)–(A3) *and is reversible with respect to its invariant probability measure* $\pi$. *Then, for all $f \in L^2(\pi)$, we have*

$$\left\| P^n f - \int f \, d\pi \right\|_{L^2} \le \rho^n \left\| f - \int f \, d\pi \right\|_{L^2},$$

*where $\rho$ is given by the formulas in Section* 2.2. *If additionally, the Markov chain is positive, then the formulas in Section* 2.3 *may be used.*

## 7. Relationship to existing results.

### 7.1. *Method of Meyn and Tweedie.* For convenience we restrict this discussion to the case when $C$ is an atom. The essence of these comments extends to the nonatomic case. Since $C$ is an atom, we can assume that $V(x) = 1$ for $x \in C$ and so (A2) is equivalent to

$$PV(x) \le \lambda V(x) + b\mathbb{1}_C(x),$$

where $b = K - \lambda$. Also we can take $\beta = P(x, C)$ for $x \in C$. Meyn and Tweedie [9] used an operator theory argument to reduce the problem to estimating the left-hand side in Proposition 4.2 at $r = 1$. If

$$\sup_{|z| \le 1} \left| \sum_{n=1}^{\infty} \left( P^n g(x) - \int g \, d\pi \right) z^n \right| \le M_1 V(x)$$

whenever $\|g\|_V \le 1$, then they can take $\rho = 1 - (M_1 + 1)^{-1}$. Using the regenerative decomposition, they obtained $M_1 \le M_2 + \zeta_C M_3$, where $M_2$ and $M_3$ can be calculated efficiently in terms of $\lambda$ and $b$, and

$$\zeta_C = \sup_{|z| \le 1} \left| 1 + \sum_{n=1}^{\infty} (u_n - u_{n-1}) z^n \right| = \sup_{|z| \le 1} |(1 - z) u(z)|,$$

where $u(z)$ is the generating function for the renewal sequence $u_n = P(X_n \in C | X_0 \in C)$. With no further information about the Markov chain, they applied a splitting technique to the forward recurrence time chain associated with the renewal sequence $u_n$ to obtain

$$\zeta_C \le \frac{32 - 8\beta^2}{\beta^3} \left( \frac{K - \lambda}{1 - \lambda} \right)^2. \tag{29}$$

We can sharpen the method of Meyn and Tweedie by putting $r = 1$ in the estimate (9) from the proof of Theorem 1.3 to get the new estimate

$$\zeta_C = \sup_{|z| = 1} |(1 - z) u(z)| = \left[ \inf_{|z| = 1} |c(z)| \right]^{-1}$$

$$\le 1 + \left( 2 \log \left( \frac{L - 1}{R - 1} \right) \right) \Big/ (\beta \log R) = 1 + \left( 2 \log \left( \frac{K - \lambda}{1 - \lambda} \right) \right) \Big/ (\beta \log \lambda^{-1}). \tag{30}$$

With more information about the Markov chain, Meyn and Tweedie obtained better estimates for $\zeta_C$. However, as they observed in [9], their method of using estimates at the value $r = 1$ to obtain estimates for $r > 1$ is very far from sharp. In particular, it cannot yield the estimate (2). By contrast, we use a version of Kendall's theorem to estimate $\sup_{|z| \le r} |(1 - z) u(z)|$ and use this together with the regenerative decomposition to estimate the left-hand side of Proposition 4.2 for $r > 1$ directly.

7.2. *Coupling method.* Our method uses (A1) and (A2) to obtain estimates on the generating function for the regeneration time $T$ for the split chain defined in (17). The estimates are based on the fact that the split chain regenerates with probability $\tilde{\beta}$ whenever $X_n \in C$. The estimate on $\mathbf{E}(r^T | X_1 \sim \nu)$, which is valid for $r < R_0$, is used with (A3) in Theorem 3.2 or 3.3 or Corollary 3.1 to obtain $\rho_C$, and then we take $\rho = \min(R_0^{-1}, \rho_C)$.

The estimates on the generating function for $T$ appear also in [15], where $R_0$ is denoted $\beta_{\mathrm{RT}}$.

The coupling method, introduced by Rosenthal [18], builds a bivariate process $\{(X_n, X'_n) : n \geq 0\}$, where each component is a copy of the original Markov chain. The stopping time of interest is the coupling time $\widehat{T} = \inf\{n \geq 0 : X_n = X'_n\}$. The minorization condition (A1) implies that the bivariate process can be constructed so that

$$P(X_{n+1} = X'_{n+1} | (X_n, X'_n) \in C \times C) \geq \tilde{\beta}.$$

Therefore, coupling can be achieved with probability $\tilde{\beta}$ whenever $(X_n, X'_n) \in C \times C$. To obtain estimates on the distribution of $\widehat{T}$, a drift condition for the bivariate process is needed. If the Markov chain is stochastically monotone and $C$ is a bottom or top set, then the univariate drift condition (A2) is sufficient. The bivariate process can be constructed so the estimates for the (univariate) regeneration time $T$ apply equally to the (bivariate) coupling time $\widehat{T}$. Thus we get $\rho = R_0^{-1}$. In particular, if $C$ is an atom, we get $\rho = \lambda$. See [7] and [21] for the case when $C$ is an atom, and [16] for the general case.

In the absence of stochastic monotonicity, a drift condition for the bivariate process can be constructed using the function $V$ which appears in (A2), but at the cost of possibly enlarging the set $C$ and also enlarging the effective value of $\lambda$. Let $b = \sup_{x \in C} PV(x) - \lambda V(x)$, so that $PV(x) \leq \lambda V(x) + b\mathbb{1}_C(x)$ for all $x \in S$. If $h(x, y) = [V(x) + V(y)]/2$, then

$$(P \times P)h(x, y) \leq \lambda_1 h(x, y) \qquad \text{if } (x, y) \notin C \times C,$$

where

$$\lambda_1 = \lambda + \frac{b}{1 + \min\{V(x) : x \notin C\}}.$$

Whereas (A2) asserts $\lambda < 1$, the coupling method requires the stronger condition $\lambda_1 < 1$. This can be achieved by enlarging the set $C$ so as to make $\min\{V(x) : x \notin C\}$ sufficiently large. Note that the condition $PV(x) \leq \lambda V(x) + b\mathbb{1}_C(x)$ for all $x \in S$ remains true with the same values of $\lambda$ and $b$ when $C$ is enlarged. However, the value of $K = \sup_{x \in C} PV(x)$ may have increased, and the value of $\tilde{\beta}$ in the minorization condition (A1) may have decreased. The coupling method now gives $\rho = \widehat{R}_0^{-1}$, where $\widehat{R}_0$ is calculated similarly to $R_0$ except that $\lambda_1$ is used in place of $\lambda$.

Here we have followed the "simple account" of the coupling method described in [19]. The assertion $\rho = \widehat{R}_0^{-1}$ is a direct consequence of [19], Theorem 1. For various developments and extensions of this method, see also [2, 4, 15].

Compared with the coupling method, our method has the advantage of allowing the use of a smaller set $C$ and a smaller numerical value of $\lambda$.

It has the disadvantage of having to apply a version of Kendall's theorem to calculate $\rho_C$. In the general setting this is a major disadvantage, but for reversible chains it is a minor disadvantage and for positive reversible chains it is no disadvantage at all.

## 8. Numerical examples.

8.1. *Reflecting random walk.* Meyn and Tweedie ([9], Section 8) considered the Bernoulli random walk on $\mathbb{Z}^+$ with transition probabilities $P(i, i - 1) = p > 1/2$, $P(i, i + 1) = q = 1 - p$ for $i \geq 1$ and boundary conditions $P(0,0) = p$, $P(0,1) = q$. Taking $C = \{0\}$ and $V(i) = (p/q)^{i/2}$, we get $\lambda = 2\sqrt{pq}$, $K = p + \sqrt{pq}$ and $\beta = p$.

For each of the values $p = 2/3$ and $p = 0.9$ considered in [9] we calculate $\rho$ in six different ways (see Table 1). Method MT is the original calculation in [9], using their formula (29) for $\zeta_C$. Method MTB is the same as MT but with our formula (30) in place of (29). Method 1.1 uses Theorem 1.1. So far these calculations have used only the values of $\lambda$, $K$ and $\beta$. The next three methods all use some extra information about the Markov chain. Method MT* uses [9] with a sharper estimate for $\zeta_C$ using the extra information that $P(\tau = 1) = p$, $P(\tau = 2) = pq$ and $\pi(0) = 1 - q/p$. Method 1.2 uses Theorem 1.2 with the extra information that the Markov chain is reversible. Finally Method LT uses the fact that the chain is stochastically monotone and gives the optimal result $\rho = \lambda$, due to Lund and Tweedie [7].

8.2. *Metropolis–Hastings algorithm for the normal distribution.* Here we consider the Markov chain that arises when the Metropolis–Hastings algorithm with candidate transition probability $q(x, \cdot) = N(x, 1)$ is used to simulate the standard normal distribution $\pi = N(0, 1)$. This example was studied by Meyn and Tweedie [9]. It also appeared in [15] and [14], where the emphasis was on convergence of the ergodic average $(1/n) \sum_{k=1}^{n} P_k(x, \cdot)$.

TABLE 1

| | $p = 2/3$ | | $p = 0.9$ | |
|---|---|---|---|---|
| | $\rho$ | $\zeta_C$ | $\rho$ | $\zeta_C$ |
| MT | 0.99994 | 1119 | 0.9967 | 78.77 |
| MTB | 0.9991 | 63.55 | 0.9470 | 2.764 |
| 1.1 | 0.9994 | | 0.9060 | |
| MT* | 0.9965 | 13 | 0.9722 | 7.313 |
| 1.2 | 0.9428 | | 0.6 | |
| LT | 0.9428 | | 0.6 | |

We compare the calculation of Meyn and Tweedie with estimates obtained by the coupling method and by our analysis. Since the Hastings–Metropolis algorithm is by construction reversible, we can use Theorem 1.2. Moreover, by Lemma 3.1 we can also apply Theorem 1.3. The continuous part of the transition probability $P(x, \cdot)$ has density

$$p(x, y) = \begin{cases} \dfrac{1}{\sqrt{2\pi}} \exp\left(-\dfrac{(y-x)^2}{2}\right), & \text{if } |x| \geq |y|, \\[2mm] \dfrac{1}{\sqrt{2\pi}} \exp\left(-\dfrac{(y-x)^2 + y^2 - x^2}{2}\right), & \text{if } |x| \leq |y|. \end{cases}$$

We use the same family of functions $V(x) = e^{s|x|}$ and sets $C = [-d, d]$ as used in [9]. Following [9] we get, for $x, s \geq 0$,

$$\lambda(x, s) := \frac{PV(x)}{V(x)}$$

$$= \exp\left(\frac{s^2}{2}\right)[\Phi(-s) - \Phi(-x - s)]$$

$$+ \exp\left(\frac{s^2}{2} - 2sx\right)[\Phi(-x + s) - \Phi(-2x + s)]$$

$$+ \frac{1}{\sqrt{2}} \exp\left(\frac{(x-s)^2}{4}\right)\Phi\left(\frac{s-x}{\sqrt{2}}\right)$$

$$+ \frac{1}{\sqrt{2}} \exp\left(\frac{x^2 - 6xs + s^2}{4}\right)\Phi\left(\frac{s-3x}{\sqrt{2}}\right)$$

$$+ \Phi(0) + \Phi(-2x) - \frac{1}{\sqrt{2}} \exp\left(\frac{x^2}{4}\right)\left[\Phi\left(\frac{-x}{\sqrt{2}}\right) + \Phi\left(\frac{-3x}{\sqrt{2}}\right)\right],$$

where $\Phi$ denotes the standard normal distribution function. Then

$$\lambda = \min_{|x| \geq d} \lambda(x, s) = \lambda(d, s), \qquad K = \max_{|x| \leq d} PV(x) = PV(d) = e^{sd}\lambda(d, s)$$

and

$$b = \max_{|x| \leq d} PV(x) - \lambda V(x) = PV(0) - \lambda V(0) = \lambda(0, s) - \lambda.$$

The computed value for $\rho$ depends on the choices of $d$ and $s$. In Table 2 we give optimal values for $d$ and $s$, and the corresponding value for $1 - \rho$ for five different methods of calculation. The first line is the calculation reported by Meyn and Tweedie, using a minorization condition with the measure $\nu$ given by

$$\nu(dx) = c \cdot \exp(-x^2)\mathbb{1}_C(x)\, dx$$

Table 2

|            | $d$  | $s$                | $1 - \rho$            |
|------------|------|--------------------|-----------------------|
| MT         | 1.4  | $4 \times 10^{-5}$ | $1.6 \times 10^{-8}$  |
| Theorem 1.1| 1    | 0.13               | $6.3 \times 10^{-7}$  |
| Coupling   | 1.8  | 1.1                | 0.00068               |
| Theorem 1.2| 1    | 0.07               | 0.0091                |
| Theorem 1.3| 1.1  | 0.16               | 0.0253                |

for a suitable normalizing constant $c$. In this case, $\nu(C) = 1$ and we have $\beta = \tilde{\beta} = \sqrt{2}\exp(-d^2)[\Phi(\sqrt{2}d) - 1/2]$. For the purposes of comparison, the other four lines were calculated using the same measure.

In Table 3, we used the measure $\nu$ given by

$$\tilde{\beta}\nu(dx) = \inf_{y \in C} p(y,x)\,dx = \begin{cases} \dfrac{1}{\sqrt{2\pi}}\exp\left(\dfrac{-(|x|+d)^2}{2}\right)dx, & \text{if } |x| \le d, \\[2mm] \dfrac{1}{\sqrt{2\pi}}e^{-d|x|-|x|^2}\,dx, & \text{if } |x| \ge d. \end{cases}$$

Now $\beta = 2[\Phi(2d) - \Phi(d)]$ and $\tilde{\beta} = \beta + \sqrt{2}\exp(d^2/4)[1 - \Phi(3d/\sqrt{2})]$. In the calculations for Theorems 1.1 and 1.2 we also used the extra information that

$$\widetilde{K} = \nu(C) + \int_{S \setminus C} V(x)\,d\nu(x) = \frac{\beta}{\tilde{\beta}} + \frac{\sqrt{2}}{\tilde{\beta}}\exp\left(\frac{(d-s)^2}{4}\right)\left[1 - \Phi\left(\frac{3d-s}{\sqrt{2}}\right)\right]$$

in the formula for $\alpha_2$.

REMARK 8.1. For this particular example, it can be verified that the process $\{|X_n| : n \ge 0\}$ is a stochastically monotone Markov chain. The coupling result of Roberts and Tweedie ([16], Theorem 2.2) can be adapted to this situation. The calculation for $\rho$ given by [16] is identical with the calculation for Theorem 1.3.

Table 3

|            | $d$  | $s$  | $1 - \rho$           |
|------------|------|------|-----------------------|
| Theorem 1.1| 1    | 0.16 | $1.7 \times 10^{-6}$  |
| Coupling   | 1.9  | 1.1  | 0.00187               |
| Theorem 1.2| 1    | 0.11 | 0.0135                |
| Theorem 1.3| 1.1  | 0.22 | 0.0333                |

8.3. *Contracting normals.* Here we consider the family of Markov chains with transition probability $P(x, \cdot) = N(\theta x, 1 - \theta^2)$ for some parameter $\theta \in (-1, 1)$. This family of examples occurs in [18] as one component of a two-component Gibbs sampler. The convergence of ergodic averages for this family was studied in [14] and [15]. Since the Markov chain is reversible with respect to its invariant probability $N(0, 1)$, we can apply Theorem 1.2. We compare these results with the estimates obtained using the coupling method.

We take $V(x) = 1 + x^2$ and $C = [-c, c]$. Then (A2) is satisfied with $\lambda = \theta^2 + 2(1 - \theta^2)/(1 + c^2)$ and $K = 2 + \theta^2(c^2 - 1)$. Also $b = \sup_{x \in C} PV(x) - \lambda V(x) = 2(1 - \theta^2)c^2/(1 + c^2)$. To ensure $\lambda < 1$, we require $c > 1$. For the minorization condition, we look for a measure $\nu$ concentrated on $C$, so that $\beta = \tilde{\beta}$. We choose $\tilde{\beta}$ and $\nu$ so that

$$\tilde{\beta}\nu(dy) = \min_{x \in C} \frac{1}{\sqrt{2\pi(1 - \theta^2)}} \exp\left(-\frac{(\theta x - y)^2}{2(1 - \theta^2)}\right) dy$$

for $y \in C$. Integrating with respect to $y$ gives

$$\tilde{\beta} = \int_{-c}^{c} \min_{x \in C} \frac{1}{\sqrt{2\pi(1 - \theta^2)}} \exp\left(-\frac{(\theta x - y)^2}{2(1 - \theta^2)}\right) dy$$

$$= 2\left[\Phi\left(\frac{(1 + |\theta|)c}{\sqrt{1 - \theta^2}}\right) - \Phi\left(\frac{|\theta|c}{\sqrt{1 - \theta^2}}\right)\right],$$

where $\Phi$ denotes the standard normal distribution function.

For the coupling method, we have $\lambda_1 = \theta^2 + 4(1 - \theta^2)/(2 + c^2)$. To ensure $\lambda_1 < 1$, we require $c > \sqrt{2}$. For the minorization condition in the coupling method there is no reason to restrict $\nu$ to be supported on $C$, so we can adapt the calculation above by integrating $y$ from $-\infty$ to $\infty$ to get

$$\tilde{\beta} = 2\left[1 - \Phi\left(\frac{|\theta|c}{\sqrt{1 - \theta^2}}\right)\right].$$

So far, the calculations have depended on $|\theta|$ but not on the sign of $\theta$. If $\theta > 0$, then $P = Q^2$, where $Q$ has parameter $\sqrt{\theta}$, so we can apply the improved estimates of Theorem 1.3. However, if $\theta < 0$, and especially if $\theta$ is close to $-1$, we can handle the almost periodicity of the chain by considering its binomial modification with transition kernel $\tilde{P} = (I + P)/2$; see [20]. Regardless of the sign of $\theta$, we can always apply Theorem 1.3 to the binomial modification. Replacing $P$ by $(1 + P)/2$ with the same $V$, $C$ and $\nu$ means replacing $\lambda$ by $(1 + \lambda)/2$, $K$ by $(1 + c^2 + K)/2$ and $\tilde{\beta}$ by $\tilde{\beta}/2$. We let $\tilde{\rho}$ denote the estimate obtained by applying Theorem 1.3 to $\tilde{P}$. Since $2n$ steps of the binomial modification $\tilde{P}$ correspond on average to $n$ steps of the original chain $P$ (see [20], Section 4), for purposes of comparison (see Table 4) we give the value of $\tilde{\rho}^2$.

| $\theta$ | Coupling | | Theorem 1.2 | | Theorem 1.3 $\theta$ positive | | Binomial mod. | |
|---|---|---|---|---|---|---|---|---|
| | $c$ | $\rho$ | $c$ | $\rho$ | $c$ | $\rho$ | $c$ | $\widehat{\rho}^2$ |
| 0.5 | 2.1 | 0.946 | 1.5 | 0.950 | 1.5 | 0.897 | 1.5 | 0.952 |
| 0.75 | 1.7 | 0.9963 | 1.2 | 0.9958 | 1.2 | 0.9847 | 1.2 | 0.9924 |
| 0.9 | 1.5 | 0.99998 | 1.1 | 0.99998 | 1.1 | 0.99948 | 1.1 | 0.99974 |

8.4. *Reflecting random walk, continued.* Here we consider the same random walk as in Section 8.1 except that the boundary transition probabilities are changed. We redefine $P(0, \{0\}) = \varepsilon$ and $P(0, \{1\}) = 1 - \varepsilon$ for some $\varepsilon > 0$. If $\varepsilon \geq p$, the Markov chain is stochastically monotone and the results of Lund and Tweedie [7] apply. Here we concentrate on the case $\varepsilon < p$, which was studied by Roberts and Tweedie [15] and Fort [4].

To apply Theorem 1.2, we take $V(i) = (p/q)^{i/2}$ and $C = \{0\}$ as earlier. Then $\lambda = 2\sqrt{pq}$, $K = \varepsilon + (1 - \varepsilon)\sqrt{p/q}$ and $\beta = \varepsilon$. If $K \leq \lambda + 2\varepsilon$ [equivalently $\varepsilon \geq (p - q)/(1 + \sqrt{q/p})$], we get $\rho = \lambda = 2\sqrt{pq}$. If $\varepsilon < (p - q)/(1 + \sqrt{q/p})$, then we take $\rho = R^{-1}$, where $R$ solves $1 + 2\varepsilon R = R^{1 + (\log K)/(\log \lambda^{-1})}$.

For the coupling method, the size of the set $C$ depends on the values of $p$ and $\varepsilon$. For the set $C = \{0, \ldots, k\}$, the condition $\lambda_1 < 1$ will be satisfied if and only if $\varepsilon > 1 - (p/q)^{k/2}(p - \sqrt{pq})$. In particular, if $\varepsilon \leq 1 - (p/q)(p - \sqrt{pq})$, then $C \supseteq \{0, 1, 2\}$ and there is no minorization condition for the time 1 transition probabilities on $C$. Instead, as pointed out in [15], it is necessary to use a minorization condition for the $m$-step kernel. This program was recently carried out by Fort. In Table 5 we denote Fort's estimates (taken from [4]) by $\rho_F$ and our estimates using Theorem 1.2 by $\rho$.

In this example, we can also calculate the exact value for $\rho_V$. We have

$$b(z) = G(z, 0) = \varepsilon z + (1 - \varepsilon) z G(z, 1)$$

$$= \varepsilon z + \frac{(1 - \varepsilon)}{2q}[1 - (1 - 4pqz^2)^{1/2}]$$

for $|z| < 1/\sqrt{4pq}$, where the formula for $G(z, 1)$ is taken from [3], Section XIV.4. The equation $b(z) = 1$ can now be solved explicitly for $|z| < 1/\sqrt{4pq}$. One solution is $z = 1$. The only other possible solution is in the interval $(-1/\sqrt{4pq}, -1)$ and exists as long as $b(-1/\sqrt{4pq}) > 1$ [equivalently as long as $\varepsilon < (p - q)/(1 + \sqrt{q/p})$]. If this condition is satisfied, the second solution is at $r = -(p - \varepsilon)/[pq + (p - \varepsilon)^2]$. By the argument in Kendall's

theorem, we deduce

$$
\rho_C = \begin{cases} \dfrac{pq + (p - \varepsilon)^2}{p - \varepsilon}, & \text{if } \varepsilon < \dfrac{p - q}{1 + \sqrt{q/p}}, \\ 2\sqrt{pq}, & \text{otherwise.} \end{cases}
$$

By inspection of this formula we see $\rho_C \geq \lambda$. Since $\rho_C \leq \rho_V \leq \max(\lambda, \rho_C)$ from (2), we deduce that $\rho_V = \rho_C$ in this example.

As $\varepsilon \to 0$, the chain becomes closer and closer to a period 2 chain. This is the setting where the binomial modification with kernel $\tilde{P} = (I + P)/2$ should converge significantly faster than the original chain: see [20]. Keeping the same function $V(x)$ and $C = \{0\}$, and applying Theorem 1.3, we get the optimal result $\tilde{\rho} = \tilde{\lambda} = (1 + \lambda)/2 = 1/2 + \sqrt{pq}$ for all $\varepsilon \geq 0$. For the purposes of comparison (see Table 6), we give the values of $\tilde{\rho}^2$ for the values of $p$ which appeared in Table 5.

## APPENDIX

PROOF OF PROPOSITION 4.1. We write $\mathcal{F}_n = \sigma\{X_r : 0 \leq r \leq n\}$. For $m \geq 0$, we have

$$
\lambda^{-1} \mathbf{E}^x (V(X_{m+1}) \mathbb{1}_{X_{m+1} \notin C} | \mathcal{F}_m) + \lambda^{-1} \mathbf{E}^x (V(X_{m+1}) \mathbb{1}_{X_{m+1} \in C} | \mathcal{F}_m) \leq V(X_m)
$$

TABLE 5

|  | $\varepsilon$ | | | $\varepsilon$ | | | $\varepsilon$ | | |
|---|---|---|---|---|---|---|---|---|---|
|  | **0.05** | **0.25** | **0.5** | **0.05** | **0.25** | **0.5** | **0.05** | **0.25** | **0.5** |
|  |  | $p = 0.6$ | |  | $p = 0.7$ | |  | $p = 0.8$ | |
| $\rho_F$ | 0.9997 | 0.9995 | 0.9994 | 0.9964 | 0.9830 | 0.9757 | 0.9793 | 0.9333 | 0.9333 |
| $\rho$ | 0.9909 | 0.9798 | 0.9798 | 0.9830 | 0.9165 | 0.9165 | 0.9759 | 0.8796 | 0.8000 |
| $\rho_V$ | 0.9864 | 0.9798 | 0.9798 | 0.9731 | 0.9165 | 0.9165 | 0.9633 | 0.8409 | 0.8000 |
|  |  | $p = 0.9$ | |  | $p = 0.95$ | |  |  |  |
| $\rho_F$ | 0.9696 | 0.8539 | 0.7500 | 0.9564 | 0.7853 | 0.5814 |  |  |  |
| $\rho$ | 0.9687 | 0.8470 | 0.6817 | 0.9645 | 0.8289 | 0.6667 |  |  |  |
| $\rho_V$ | 0.9559 | 0.7885 | 0.6250 | 0.9528 | 0.7679 | 0.5556 |  |  |  |

TABLE 6

| $p$ | **0.6** | **0.7** | **0.8** | **0.9** | **0.95** |
|---|---|---|---|---|---|
| $\tilde{\rho}^2$ | 0.9799 | 0.9186 | 0.8100 | 0.6400 | 0.5154 |

on the set $\{X_m \notin C\}$. Multiply by $\lambda^{-m}\mathbb{1}_{\tau>m}$, take expectation and sum over $m = 0$ to $n - 1$ to obtain

$$(31) \qquad \lambda^{-n}\mathbf{E}^x(V(X_n)\mathbb{1}_{\tau>n}) + \mathbf{E}^x(\lambda^{-\tau}V(X_\tau)\mathbb{1}_{\tau\le n}) \le V(x) \qquad \text{for all } x \notin C$$

or, equivalently,

$$\lambda^{-n}\mathbf{E}^x(V(X_n)\mathbb{1}_{\tau\ge n}) + \mathbf{E}^x(\lambda^{-\tau}V(X_\tau)\mathbb{1}_{\tau<n}) \le V(x) \qquad \text{for all } x \notin C.$$
(32)

This implies that $\mathbf{P}^x(\tau \ge n) \le \lambda^n V(x)$ for $x \notin C$, which implies (i).

The first assertion in (ii) is obtained by letting $n \to \infty$ in (31) and the second assertion follows from the first via the identity

$$G(r,x) = rP(x,C) + r\int_{S\setminus C} P(x,dy)G(r,y).$$

For the calculations to prove (iii) and (iv) it is convenient to define the function

$$J(r,x) = \mathbf{E}^x(r^\tau V(X_\tau)).$$

The functions $H$ and $J$ satisfy the identities

$$(33) \qquad H(r,x) = rPV(x) + r\int_{S\setminus C} P(x,dy)H(r,y)$$

and

$$(34) \qquad J(r,x) = r\int_C P(x,dy)V(y) + r\int_{S\setminus C} P(x,dy)J(r,y).$$

For $0 < r < \lambda^{-1}$, multiply (32) by $\lambda^n r^n$ and sum over $n = 1$ to $\infty$. We obtain

$$(35) \qquad H(r,x) + \frac{\lambda r}{1 - \lambda r}J(r,x) \le \frac{\lambda r}{1 - \lambda r}V(x) \qquad \text{for all } x \notin C,$$

which gives the first part of (iii). For $x \in C$, we use the inequality (35) in the right-hand side of the identity (33) along with the identity (34) to obtain

$$H(r,x) \le rPV(x) + \frac{\lambda r^2}{1 - \lambda r}\int_{S\setminus C} P(x,dy)[V(y) - J(r,y)]$$

$$= \frac{r}{1 - \lambda r}PV(x) - \frac{\lambda r}{1 - \lambda r}J(r,x)$$

$$\le \frac{r(K - \lambda r)}{1 - \lambda r}.$$

This completes (iii). If we replace $\lambda^n r^n$ by $\lambda^n(r^n - 1)$ in the derivation of (35) we obtain instead

$$H(r,x) - H(1,x) + \frac{\lambda r}{1 - \lambda r}J(r,x) - \frac{\lambda}{1 - \lambda}J(1,x)$$
(36)
$$\le \frac{\lambda(r-1)}{(1-\lambda)(1-r\lambda)}V(x)$$

for $x \notin C$ and $1 < r < \lambda^{-1}$. Using (33), (36) and (34), we get

$$H(r,x) - H(1,x)$$

$$= r \int_{S \setminus C} P(x,dy)[H(r,y) - H(1,y)]$$

$$\leq \frac{\lambda r(r-1)}{(1-\lambda)(1-\lambda r)} \int_{S \setminus C} P(x,dy) V(y)$$

$$- r \int_{S \setminus C} P(x,dy) \left[ \frac{\lambda r}{1 - \lambda r} J(r,y) + \frac{\lambda}{1-\lambda} J(1,y) \right]$$

$$= \frac{\lambda r(r-1)}{(1-\lambda)(1-r\lambda)} PV(x) - \lambda r \left[ \frac{1}{1-\lambda r} J(r,x) - \frac{1}{1-\lambda} J(1,x) \right]$$

and (iv) follows easily. $\quad\square$

PROOF OF PROPOSITION 4.2. For $z \in \mathbf{C}$, write

$$G(z,x) = \mathbf{E}^x(z^\tau), \qquad H(z,x) = \mathbf{E}^x \left( \sum_{n=1}^{\tau} z^n V(X_n) \right)$$

and

$$H_g(z,x) = \mathbf{E}^x \left( \sum_{n=1}^{\tau} z^n g(X_n) \right).$$

Let $u(z) = \sum_{n=0}^{\infty} u_n z^n$ be the generating function for the sequence $u_n$. Suppose $|z| < 1$. The first-entrance–last-exit decomposition ([8], equation (13.46)) yields

$$\sum_{n=1}^{\infty} P^n g(x) z^n = H_g(z,x) + G(z,x) u(z) H_g(z,a).$$

Furthermore, [8], equation (13.50), gives

$$\int g \, d\pi = \pi(C) H_g(1,a).$$

Together, for $|z| < 1$ we have

$$\sum_{n=1}^{\infty} \left( P^n g(x) - \int g \, d\pi \right) z^n$$

$$= H_g(z,x) + G(z,x) u(z) H_g(z,a) - \frac{z\pi(C)}{1-z} H_g(1,a)$$

(37)

$$= H_g(z,x) + G(z,x) \left[ u(z) - \frac{\pi(C)}{1-z} \right] H_g(z,a)$$

$$- \pi(C) H_g(z,a) \frac{G(z,x) - 1}{z-1} - \pi(C) \frac{H_g(z,a) - z H_g(1,a)}{z-1}.$$

Now

$$\left| \frac{H_g(z,a) - zH_g(1,a)}{z-1} \right|$$

$$= \left| \mathbf{E}^a\left( \sum_{n=1}^{\tau} g(X_n)(z + \cdots + z^{n-1}) \right) \right|$$

$$\leq \mathbf{E}^a\left( \sum_{n=1}^{\tau} V(X_n)(|z| + \cdots + |z|^{n-1}) \right)$$

$$\leq \frac{H(r,a) - rH(1,a)}{r-1}$$

if $|z| \leq r$ and $r > 1$, and a similar estimate holds for $|(G(z,x) - 1)/(z-1)|$. Also $\pi(C) = \lim_{n\to\infty} \mathbf{P}^a(X_n \in C) = \lim_{n\to\infty} u_n = u_\infty$, so

$$u(z) - \frac{\pi(C)}{1-z} = \sum_{n=0}^{\infty} (u_n - u_\infty)z^n$$

and the result now follows easily from (37). $\square$

PROOF OF PROPOSITION 4.3. Notice that the invariant probability measure $\pi$ for $\{X_n : n \geq 0\}$ is the $S$ marginal of the stationary probability $\bar{\pi}$, say, for the split chain, so that $\int g\, d\bar{\pi} = \int g\, d\pi$. The argument used in the proof of Proposition 4.2 gives expressions similar to (37) for $\sum_{n=1}^{\infty}(\overline{\mathbf{E}}^{x,i}(g(X_n)) - \int g\, d\pi)z^n$ for $i = 0, 1$. Multiplying the $i = 0$ expression by $(1 - \tilde{\beta})\mathbb{1}_C(x)$ and the $i = 1$ expression by $\tilde{\beta}\mathbb{1}_C(x)$ and adding gives an expression for $\sum_{n=1}^{\infty}(P^n g(x) - \int g\, d\pi)z^n$. The remainder of the proof is exactly as in the proof of Proposition 4.2. $\square$

To prove Proposition 4.4 we need some intermediate results. Define

$$G(r,x,i) = \overline{\mathbf{E}}^{x,i}(r^\tau),$$

$$H(r,x,i) = \overline{\mathbf{E}}^{x,i}\left( \sum_{n=1}^{\tau} r^n V(X_n) \right).$$

In addition to $\widetilde{G}(r)$ defined in Section 4.2, we define

$$\widetilde{H}(r) = \sup\{H(r,x,0) : x \in C\},$$

$$\widetilde{H}(r,1) = \sup\{H(r,x,0) - rH(1,x,0) : x \in C\}.$$

We need to consider the following functions which are defined in terms of the split chain and the original stopping time $\tau = \inf\{n \geq 1 : X_n \in C\}$.

LEMMA A.1. *Assume conditions* (A1) *and* (A2). *Then*

$$(38) \qquad \overline{G}(r,x,i) \leq \frac{\tilde{\beta}G(r,x,i)}{1-(1-\tilde{\beta})\widetilde{G}(r)},$$

$$(39) \qquad \overline{H}(r,x,i) \leq H(r,x,i) + \frac{(1-\tilde{\beta})\widetilde{H}(r)G(r,x,i)}{1-(1-\tilde{\beta})\widetilde{G}(r)}$$

*and*

$$\overline{H}(r,x,i) - r\overline{H}(1,x,i)$$

$$(40) \qquad \leq H(r,x,i) - rH(1,x,i) + \frac{(1-\tilde{\beta})\widetilde{H}(r,1)G(r,x,i)}{1-(1-\tilde{\beta})\widetilde{G}(r)}$$

$$\qquad + \frac{(1-\tilde{\beta})r\widetilde{H}(1)}{1-(1-\tilde{\beta})\widetilde{G}(r)}\Big([G(r,x,i)-1] + \frac{(1-\tilde{\beta})}{\tilde{\beta}}[\widetilde{G}(r)-1]\Big)$$

*for all* $r > 1$ *such that* $(1-\tilde{\beta})\widetilde{G}(r) < 1$ *and* $r < \lambda^{-1}$.

PROOF. Define the sequence of stopping times $\tau_0 = 0$ and $\tau_k = \tau_{k-1} + \tau \circ \theta(\tau_{k-1})$ for $k \geq 1$ [where $\theta(n)$ denotes the natural time $n$ shift]. Define the random variable $K = \inf\{k \geq 1 : Y_{\tau_k} = 1\}$, so that $T = \tau_K$. Then

$$\overline{\mathbf{E}}^{x,i}\bigg(\sum_{n=1}^{T} r^n V(X_n)\bigg)$$

$$= \overline{\mathbf{E}}^{x,i}\bigg(\sum_{k=1}^{K} \sum_{n=\tau_{k-1}+1}^{\tau_k} r^n V(X_n)\bigg)$$

$$= \sum_{k=1}^{\infty} \overline{\mathbf{E}}^{x,i}\bigg(\sum_{n=\tau_{k-1}+1}^{\tau_k} r^n V(X_n), K \geq k\bigg)$$

$$= H_r(x,i) + \sum_{k=2}^{\infty} \overline{\mathbf{E}}^{x,i}\bigg(\sum_{n=\tau_{k-1}+1}^{\tau_k} r^n V(X_n), K \geq k\bigg).$$

By conditioning on $\mathcal{G}(\tau_{k-1})$, where $\mathcal{G}(n) = \mathcal{F}_n^X \vee \mathcal{F}_{n-1}^Y$, we get

$$\overline{\mathbf{E}}^{x,i}\bigg(\sum_{n=\tau_{k-1}+1}^{\tau_k} r^n V(X_n), K \geq k\bigg) \leq (1-\tilde{\beta})\widetilde{H}(r)\overline{\mathbf{E}}^{x,i}(r^{\tau_{k-1}}, K \geq k-1)$$

and

$$\overline{\mathbf{E}}^{x,i}(r^{\tau_k}, K \geq k) \leq (1-\tilde{\beta})\widetilde{G}(r)\overline{\mathbf{E}}^{x,i}(r^{\tau_{k-1}}, K \geq k-1)$$

for $k \geq 2$. Together we obtain by induction

$$\overline{\mathbf{E}}^{x,i}\left(\sum_{n=1}^{T} r^n V(X_n)\right)$$

$$\leq H_r(x,i) + (1-\tilde{\beta})\tilde{H}(r)\sum_{k=1}^{\infty}\overline{\mathbf{E}}^{x,i}(r^{\tau_k}, K \geq k)$$

$$\leq H_r(x,i) + \tilde{H}(r)\sum_{k=1}^{\infty}(1-\tilde{\beta})^k \tilde{G}(r)^{k-1}G(r,x,i),$$

giving (39). To prove (38), note first that

$$\overline{\mathbf{E}}^{x,i}(r^T) = \sum_{k=1}^{\infty}\overline{\mathbf{E}}^{x,i}(r^{\tau_k}, K=k) = \tilde{\beta}\sum_{k=1}^{\infty}\overline{\mathbf{E}}^{x,i}(r^{\tau_k}, K \geq k),$$

and the remainder of the proof is the special case of the proof above with $V$ replaced by $\mathbb{1}_C$. To prove (40), we note that for $k \geq 2$,

$$\overline{\mathbf{E}}^{x,i}\left(\sum_{n=\tau_{k-1}+1}^{\tau_k}(r^n-r)V(X_n), K \geq k\right)$$

$$\leq (1-\tilde{\beta})\tilde{H}(r,1)\overline{\mathbf{E}}^{x,i}(r^{\tau_{k-1}}, K \geq k-1)$$

$$+ (1-\tilde{\beta})\tilde{H}(1)r\overline{\mathbf{E}}^{x,i}(r^{\tau_{k-1}}-1, K \geq k-1)$$

and $\overline{\mathbf{P}}^{x,i}(K \geq k-1) = (1-\tilde{\beta})^{k-2}$. Then the rest of the proof is essentially the same as for (39). □

LEMMA A.2. *Assume conditions* (A1) *and* (A2), *and let* $\alpha_1$ *and* $\alpha_2$ *be given by* (18) *and* (20) *of Proposition* 4.4. *Then for* $1 < r \leq \lambda^{-1}$,

$$\tilde{G}(r) \leq r^{\alpha_1},$$

$$G(r,a,1) \leq r^{\alpha_2}.$$

PROOF.   For $x \in C$ we have

$$(1-\tilde{\beta})G(\lambda^{-1}, x, 0)$$

$$= \lambda^{-1}\left[P(x,C) - \tilde{\beta}\nu(C) + \int_{S\setminus C}G(\lambda^{-1},y)[P(x,dy) - \tilde{\beta}\nu(dy)]\right]$$

$$\leq \lambda^{-1}\left[P(x,C) - \tilde{\beta}\nu(C) + \int_{S\setminus C}V(y)[P(x,dy) - \tilde{\beta}\nu(dy)]\right]$$

$$\leq \lambda^{-1}\left[PV(x) - \tilde{\beta}\int_S V\,d\nu\right]$$

$$\leq \lambda^{-1}(K - \tilde{\beta}).$$

Therefore, for $1 < r \leq \lambda^{-1}$,

$$
\begin{aligned}
\widetilde{G}(r) = \sup_{x \in C} G(r, x, 0) &\leq \sup_{x \in C} (G(\lambda^{-1}, x, 0))^{(\log r)/(\log \lambda^{-1})} \\
&\leq \left( \frac{\lambda^{-1}(K - \tilde{\beta})}{1 - \tilde{\beta}} \right)^{(\log r)/(\log \lambda^{-1})} \\
&= r^{\alpha_1}.
\end{aligned}
$$

(This estimate on $\widetilde{G}(r)$ appears as Theorem 2.2 in [15].) The minorization condition implies $\tilde{\beta} \int_S V \, d\nu \leq PV(x) \leq K$ for $x \in C$ and so $\int_S V \, d\nu \leq K/\tilde{\beta}$. We have

$$
\begin{aligned}
G(\lambda^{-1}, a, 1) &= \lambda^{-1} \nu(C) + \lambda^{-1} \int_{S \setminus C} G(\lambda^{-1}, y) \nu(dy) \\
&\leq \lambda^{-1} \int_S V(y) \nu(dy) \\
&\leq \frac{K}{\lambda \tilde{\beta}}
\end{aligned}
$$

and so, for $1 < r \leq \lambda^{-1}$,

$$
G(r, a, 1) \leq \left( \frac{K}{\lambda \tilde{\beta}} \right)^{(\log r)/(\log \lambda^{-1})} = r^{\alpha_2}
$$

and the proof is complete. $\quad\square$

PROOF OF PROPOSITION 4.4 AND REMARK 4.1. It is clear from the proof of Lemma A.2 that its assertions remain valid when $\alpha_2$ is chosen according to Remark 4.1. The inequality (19) is part of the statement of Lemma A.2, and inequalities (21) and (22) are immediate consequences of Lemmas A.1 and A.2. The result (23) uses the estimate

$$
(1 - \tilde{\beta}) \widetilde{H}(r) \leq \sup_{x \in C} H(r, x) - \tilde{\beta} r \leq \frac{r[K - r\lambda - \tilde{\beta}(1 - r\lambda)]}{1 - r\lambda}
$$

from Lemma A.1. To obtain (24), notice first that

$$
\begin{aligned}
&H(r, a, 1) + \frac{(1 - \tilde{\beta}) \widetilde{H}(r) G(r, a, 1)}{1 - (1 - \tilde{\beta}) \widetilde{G}(r)} \\
&= \frac{1}{1 - (1 - \tilde{\beta}) \widetilde{G}(r)} \left[ G(r, a, 1) \sup_{x \in C} H(r, x) + H(r, a, 1) \left[ 1 - \sup_{x \in C} G(r, x) \right] \right] \\
&\leq \frac{1}{1 - (1 - \tilde{\beta}) \widetilde{G}(r)} G(r, a, 1) \sup_{x \in C} H(r, x).
\end{aligned}
$$

The proof of (25) is similar, using the inequality

$$H(r,a,1) - rH(1,a,1) + \frac{(1-\tilde{\beta})\widetilde{H}(r,1)G(r,a,1)}{1-(1-\tilde{\beta})\widetilde{G}(r)}$$

$$\leq \frac{1}{1-(1-\tilde{\beta})\widetilde{G}(r)}G(r,a,1)\sup_{x\in C}[H(r,y)-rH(1,y)]. \qquad \square$$

**Acknowledgments.** I am grateful to Sean Meyn and Richard Tweedie for valuable comments on a much earlier version of this paper [1], and to Gersende Fort for providing access to her preprint [4]. In particular, I wish to thank Sean Meyn for pointing out the result in Corollary 3.1.

DEPARTMENT OF MATHEMATICS
UNIVERSITY OF SOUTHERN CALIFORNIA
LOS ANGELES, CALIFORNIA 90089-2532
USA
E-MAIL: baxendal@math.usc.edu
URL: http://math.usc.edu/˜baxendal/